\documentclass[11pt]{article}
\usepackage[latin1]{inputenc} 
\usepackage{latexsym,amsmath,amsfonts,euscript} 
\newtheorem{Th}{Theorem} 
\newtheorem{Theorem}{Theorem}[section] 
\newtheorem{Prop}[Theorem]{Proposition} 
\newtheorem{Def}[Theorem]{Definition} 
\newtheorem{Lemma}[Theorem]{Lemma} 
 
\newtheorem{Remark}[Theorem]{Remark} 
 
\begin{document}

\newcommand{\finishproof}{\hfill $\Box$ \vspace{5mm}} 
\newcommand{\cI}{{\cal I}} 
\newcommand{\al}{\alpha} 
\newcommand{\be}{\beta} 
\newcommand{\ga}{\gamma} 
\newcommand{\Ga}{\Gamma} 
\newcommand{\ep}{\epsilon} 
\newcommand{\de}{\delta} 
\newcommand{\De}{\Delta} 
\newcommand{\ka}{\kappa} 
\newcommand{\la}{\lambda}
\newcommand{\La}{\Lambda}
\newcommand{\te}{\theta} 
\newcommand{\om}{\omega} 
\newcommand{\si}{\sigma} 
\newcommand{\bg}{{\bar g}} 
\newcommand{\bs}{{\bar s}} 
\newcommand{\bL}{{\bar L}} 
\newcommand{\tL}{{\tilde L}}
\newcommand{\tK}{{\tilde K}}
\newcommand{\tQ}{{\tilde Q}}
\newcommand{\tg}{\tilde g} 
\newcommand{\eqdef}{:=} 
\newcommand{\const}{\mathop{\rm const}} 
\newcommand{\sg}{\mathop{\rm sgrad}} 
\newcommand{\id}{{\bf 1}} 
\newcommand{\dds}{\frac{d}{ds}|_{s=0}} 
\newcommand{\p}{{\partial}} 
\newcommand{\C}{\mathbb C}
\newcommand{\Z}{\mathbb Z} 
\newcommand{\N}{\mathbb N} 
\newcommand{\R}{\mathbb R}
\newcommand{\D}{\mathbb D} 
\newcommand{\T}{\mathbb T}
\newcommand{\tI}{\tilde I}
\newcommand{\s}{\mathbb S}

\renewcommand{\mod}{\mathop{\rm mod}} 
\renewcommand{\theequation}{\thesection .\arabic{equation}} 
\renewcommand{\arraystretch}{1.3}

\title{On the Integral  Geometry of Liouville Billiard Tables} 
\author{G. Popov and P. Topalov} 
\date{
} 
\maketitle

\begin{abstract} 
\noindent The notion of a Radon transform is introduced for completely integrable billiard tables.
In the case of Liouville billiard tables of dimension $3$ we prove that the Radon transform is one-to-one on
the space of continuous functions $K$ on the boundary which are  invariant with respect to the corresponding group
of symmetries. We prove also that the frequency map associated with a class of Liouville billiard tables is non-degenerate.
This allows us to obtain  spectral rigidity of  the corresponding Laplace-Beltrami operator with Robin boundary conditions.
\end{abstract} 
 
\setcounter{section}{0}
\section{Introduction} 
\setcounter{equation}{0}
This paper is concerned with the integral geometry and the spectral rigidity of Liouville billiard tables. 
By a  billiard table we mean a smooth compact connected Riemannian manifold $(X,g)$ of dimension $n \ge 2$
with a non-empty boundary $\Gamma := \partial X$.  The elastic reflection of geodesics at $\Gamma$ determines
continuous curves on $X$ called {\em billiard trajectories} as well as a discontinuous dynamical system on $T^*X$ --
the ``{\em billiard flow}'' -- that generalizes the geodesic flow on closed manifolds without boundary.
The billiard flow on $T^*X$ induces a discrete dynamical system in the open coball bundle $B^\ast\Gamma$ of $\Gamma$
given by the corresponding billiard ball map $B$  and its iterates. The  map $B$ is defined in an open  subset of 
$B^\ast\Gamma=\{\xi\in T^*\Gamma\,:\,\|\xi\|_g<1\}$, where $\|\xi\|_g$ denotes the norm induced by the Riemannian metric $g$
on the corresponding cotangent plane and it can be considered as a discrete Lagrangian systems  as in  \cite{Mos-Ves},
\cite{PT2}, \cite{Tabach}. 
The  orbits of $B$ can be obtained  by a variational principal  and they can be viewed  as ``discrete geodesics'' of
the corresponding Lagrangian. In this context, periodic orbits of $B$  can be considered as ``discrete  closed  geodesics''. 

Let $\mu$ be a positive continuous function on $B^\ast\Gamma$. Denote by $\pi_\Gamma^\ast K$ the pull-back of the continuous
function  $K\in C(\Gamma)$ with respect to the projection $\pi_\Gamma : T^*\Gamma\to\Gamma$. 
We are interested in the following problems.

\vspace{0.3cm}

\noindent{\large\bf Problem A.} Let $K$ be a continuous function on $\Gamma$ such that the mean value of the product 
$\pi_\Gamma^\ast K\cdot\mu$  is zero on any periodic orbit of the billiard ball map $B$. Does it imply $K\equiv 0$?

\vspace{0.3cm}
The mapping assigning to any periodic orbit $\gamma=\{\varrho_0,\varrho_1,\ldots,\varrho_{m-1}\}\subset B^\ast\Gamma$ of
the map $B$ the mean value $(1/m)\sum_{j=0}^{m-1}\left(\pi_\Gamma^\ast K\cdot\mu\right)(\varrho_j)$ of the function
$\pi_\Gamma^\ast K\cdot\mu$ on  $\gamma$ can be viewed as a discrete analogue of the Radon transform,  
considering  the periodic orbits of the billiard ball map as  discrete  closed geodesics. 
Problem A has a positive answer for any ball in the Euclidean space $\R^n$ centered at the origin if
$\mu=1$ and $K$ is even. In fact, approximating the great circles on the sphere  by closed billiard trajectories of the billiard table
we obtain from the hypothesis in Problem A that the integral of $K$ over any great circle is zero.
Since $K$ is even, by Funk's theorem we obtain  $K\equiv 0$ (\cite[Theorem 4.53]{Besse}).
The case of general Riemannian manifold is much more complicated. 

Denote  by $\pi_X : T^* X\to X$ the natural projection of the cotangent bundle $T^* X$ onto $X$.
Let $S^\ast X|_\Gamma=\{\xi\in T^*X\,:\,\pi_X(\xi)\in\Gamma, \|\xi\|_g=1\}$ be the restriction of the
unit co-sphere bundle to $\Gamma$.
There are two natural choices for the function $\mu$ we are concerned with, namely, $\mu \equiv 1$ or 
$\mu(\xi) = \langle \pi^+(\xi),n_g\rangle^{-1}$, $\xi\in B^\ast\Gamma$,
where $\langle\cdot,\cdot\rangle$ is the standard pairing
between vectors and covectors, $n_g$ is the inward unit normal to $\Gamma$ at $x=\pi_\Gamma(\xi)$, and 
$\pi^+ : B^\ast\Gamma\to S^\ast X|_\Gamma$ assigns to any $\xi\in T^\ast_x\Gamma$ with norm $\|\xi\|_g < 1$ 
the unit outgoing covector the restriction of which  to $T_x\Gamma$ coincides with $\xi$.
Recall that a covector based on $x$ is outgoing if its value on $n_g(x)$ is non-negative. 
The latter choice of $\mu$ is related with the wave-trace formula for manifolds with boundary
obtained by V. Guillemin and R. Melrose \cite{GM, GM1}. 
It appears also in the iso-spectral invariants of the Robin boundary problem for the Laplace-Beltrami operator
obtained in \cite{PT3}. From now on we fix the positive function  $\mu\in C(B^\ast\Gamma)$ by
\begin{equation}\label{thedensity}
 \mu\equiv 1\, ,    \quad \mbox{or by} \quad  
\mu(\xi)=\langle \pi^+(\xi),n_g\rangle^{-1}\, ,\  \xi\in B^\ast\Gamma \, .
\end{equation}
For that choice of $\mu$, it will be shown that Problem A has a positive solution for a class of  Liouville billiard tables of classical type.  
A Liouville billiard table (shortly L.B.T.) of dimension $n\ge 2$,  is a completely integrable billiard table $(X,g)$ (the notion of complete integrability will be
recalled in Sect.\,\ref{sec:set-up}) admitting $n$ functionally independent and Poisson commuting integrals of the billiard flow on $T^*X$
which are quadratic forms in the momentum. A L.B.T. can be viewed as a $2^{n-1}$-folded branched covering of a disk-like
domain in $\R^n$ by the cylinder $ \T^{n-1}\times [-N,N]$, where $\T=\R/\Z$ and $N>0$. 
Liouville billiard tables of dimension two are defined in \cite{PT1} and in any dimension $n\ge 2$ in \cite{PT2}, 
where the integrability of the billiard ball map is shown via the geodesic equivalence principal. Here we  write explicitly first integrals of the billiard flow and show that it is completely integrable (see Sect.\, \ref{sec:construction}). 
An important subclass of L.B.T.s  are the Liouville billiard tables of classical type having an  additional symmetry and for which the boundary is strictly geodesically convex (with respect to the outward normal $-n_g$). 
It turns out that the  group of isometries of a L.B.T.  of classical type is isomorphic to $(\Z/2\Z)^n$. 
Moreover, the group of isometries of $(X,g)$ induces a group of isometries $G$ on $\Gamma$ which is isomorphic to
$(\Z/2\Z)^n$. An important example of a L.B.T. of classical type is the interior of the
$n$-axial ellipsoid equipped with the Euclidean metric. More generally, there is a non-trivial  two-parameter family
of L.B.T.s of classical type  of constant scalar curvature $\kappa$ having the same broken geodesics
(considered as non-parameterized curves) as the ellipsoid \cite[Theorem 3]{PT2}.
This family includes the ellipsoid ($\kappa=0$), a L.B.T. on the sphere ($\kappa=1$) and a L.B.T.
in the hyperbolic space ($\kappa=-1$).

\begin{Th}\hspace{-2mm}{\bf .}\label{Th:Problem_A}
Let $(X,g)$, $\mbox{dim}\, X=3$, be an analytic L.B.T. of classical  type. Suppose that  there is at least one non-periodic
geodesic on the boundary $\Gamma$. Choose $\mu$ as in (\ref{thedensity}). Let  $K\in C(\Gamma)$ be invariant with
respect to the group of isometries $G \cong (\Z/2\Z)^3$ of the boundary $\Gamma$ and such  that the mean value of
$\pi_\Gamma^\ast K\cdot\mu$ on any periodic orbit of the billiard ball map is zero. Then $K\equiv 0$.  
\end{Th}
In particular, Problem $A$ has a positive solution for ellipsoidal billiard tables in $\R^3$  with $\mu\equiv 1$ as well as for 
$\mu(\xi)=\langle \pi^+(\xi),n_g\rangle^{-1}$, for any $K\in C(\Gamma)$ which is invariant under the reflections with respect to the coordinate planes
$O_{xy}$, $O_{yz}$, and $O_{xz}$. More generally, Theorem 3 can be applied for any L.B.T.  of the family described in  \cite[Theorem 3]{PT2}. The condition that the boundary contains at least one non-closed geodesic will become clear after the discussion of Problem C. 

As it was mentioned above the map assigning to each periodic orbit of the billiard ball map $B$ the mean value of
$\pi_\Gamma^\ast K\cdot\mu$ on it can be considered as a discrete analogue of the Radon transform. Another version of the Radon transform can be defined as follows. Denote by ${\mathcal F}$ the family of all  Lagrangian   tori $\Lambda\subset B^\ast\Gamma$ which are invariant with respect to  some exponent $B^m$, $m\ge 1$,  of the billiard ball map $B$, i.e. 
$B^m(\Lambda)\subseteq \Lambda$.  For any continuous function $K$ on $\Gamma$  we denote by ${\mathcal R}_{K,\mu}(\Lambda)$ the mean value of the integral
of $\pi_\Gamma^\ast K\cdot\mu$ on $\Lambda\in {\mathcal F}$ with respect to the Leray form (see Sect.\,\ref{sec:set-up}).
The mapping $\Lambda \mapsto {\mathcal R}_{K,\mu}(\Lambda)$, $\Lambda\in {\mathcal F}$, 
will be called a Radon transform of $K$ as well.

\vspace{0.3cm}

\noindent{\large\bf Problem B.} Let $K$ be a continuous function on $\Gamma$ which is invariant with respect to the group of isometries $G$.  Does the relation ${\mathcal R}_{K,\mu}\equiv 0$  imply $K\equiv 0$?

\vspace{0.3cm}

\noindent The main result of the paper is the following theorem, which gives a positive answer of Problem B for L.B.T.s.  
\begin{Th}\hspace{-2mm}{\bf .}\label{Th:Problem_B}
Let $(X,g)$, $\mbox{dim}\, X=3$, be a   Liouville billiard table of classical  type. Fix $\mu$ by (\ref{thedensity}). 
If $K\in C(\Gamma)$ is invariant under the group of symmetries $G$ of $\Gamma$ and ${\mathcal R}_{K,\mu}(\Lambda) = 0$ for any  $\Lambda\in {\mathcal F}$, then  $K\equiv 0$. 
\end{Th}
We point out that L.B.T.s of classical type are smooth by construction but they are not supposed to be analytic. 

A similar result has been obtained for the ellipse in \cite{GM} and more generally for  L.B.T.s of classical type
in dimension $n=2$ in \cite{PT1} and \cite{PT3}. 
It is always interesting to find a smaller set of data $\Lambda$ for which the Radon transform is one-to-one. 
In the case $n=2$ the proof is done by  analyticity, and we need to know the values of the Radon transform
${\mathcal R}_{K,\mu}(\Lambda)$ only on a family of invariant circles $\{\Lambda_j\}_{j\in\N}$ approaching
the boundary $S^\ast\Gamma$ of $B^\ast\Gamma$. The case $n=3$ is more complicated, since the argument using
analyticity does not work any more. Nevertheless, we can restrict the Radon transform to data ``close'' to
the boundary in the following sense:
It will be shown in Sect. \ref{subsec:parameter} that any  L.B.T. of classical type of dimension 3 admits
four not necessarily connected charts $U_j$, $1\le j\le 4$, of action-angle variables in $B^\ast\Gamma$.
Two of them, say $U_1$ and $U_2$, have the property that any unparameterized geodesic in $S^\ast\Gamma$ can be
obtained as a limit of orbits of $B$ lying either in $U_1$ or in $U_2$
(then the corresponding broken geodesics approximate geodesics of the boundary). Moreover, in any connected component of $U_1$ and $U_2$ there is such a sequence of orbits of $B$, while $U_3$ and $U_4$ do not enjoy this property. In other words,   the charts  $U_1$ and $U_2$  can be  characterized by the property that there is a family of {\em ``whispering gallery rays''}  issuing  from any of their connected components.  For this reason the  two cases $j=1,2$ will be  referred  as to boundary cases. 
Denote by ${\mathcal F}_b$ the set of all $\Lambda\in {\mathcal F}$ lying either in $U_1$ or in $U_2$. 
We will show in 
Theorem \ref{Th:R-rigid} 
that the restriction of the Radon transform ${\mathcal R}_{K,\mu}$ on ${\mathcal F}_b$ determines uniquely $K$. 

As an application we prove spectral rigidity of the Robin boundary problem for Liouville billiard tables.
Given a real-valued function $K\in C(\Gamma,\R)$, we consider the ``positive'' Laplace-Beltrami operator $\Delta$
on $X$ with domain     
\[  
D:= \left\{u\in H^2(X)\,:\,\,\,\frac{\partial u}{\partial n_g}|_\Gamma = K  
u|_\Gamma\right\}\, ,  
\]   
where $H^2(X)$ is the Sobolev space, and  $n_g(x)$, $x\in \Gamma$, is the inward unit normal to  
$\Gamma$ with respect to the metric $g$. We denote this operator by  
$\Delta_{g,K}$.   
It is a selfadjoint operator in $L^2(X)$ with  
discrete spectrum   
\[  
{\rm Spec}\,  \Delta_{g,K}:=\{ \lambda_1 \le \lambda_2 \le \cdots  
\}\,,  
\]  
where each eigenvalue $\lambda=\lambda_j$  is repeated according to its  
multiplicity, and it solves the spectral problem  
\begin{equation}\label{thespectrum}    
\left\{  
\begin{array}{rcll}  
\Delta\,  u\ &=& \ \la\,  u \, \quad \mbox{in}\ $X$\, , \\  
\displaystyle\frac{\partial u}{\partial n_g}|_\Gamma \  &=& \ K\, u|_\Gamma \, .  
\end{array}  
\right.
\end{equation}\label{isospectral}
Let  $[0,1]\ni t\mapsto K_t\in C^\infty(\Gamma,\R)$  be a continuous family of   smooth 
real-valued functions on $\Gamma$. To simplify the notations we denote by 
$\Delta_t$ the corresponding operators $\Delta_{g,K_t}$. This family is said
to be isospectral if
\begin{equation}\label{e:isospectral}
\forall\,  t \in [0,1]\, ,\  {\rm Spec}\left(\Delta_t \right)\, = \, 
{\rm Spec}\left(\Delta_0 \right)\, .
\end{equation}
We  consider here a weaker notion of isospectrality which has been introduced in \cite{PT3}. 
Fix  two positive constants  $c$ and $d>1/2$,
and  consider the union of infinitely many disjoint intervals 
\begin{enumerate}
\item[(H$_1$)]
{\em $\mathcal{I}\ :=\ \displaystyle\cup_{k=1}^{\infty}\ [a_k,b_k]$,
\quad 
$0< a_1<b_1< \cdots <a_k<b_k<\cdots $ , \quad  such that \\[0.3cm]
$\displaystyle{\lim_{k\to\infty}a_k\, =\, \lim_{k\to\infty}b_k\,  =\,  +\infty\, ,\   
\lim_{k\to\infty}(b_k -a_k)\,  =\,  0\, ,\  
\mbox{and} \ a_{k+1} - b_{k}\, \ge \, c b_k^{-d}\ \mbox{for  any} \
k\ge 1\, .}$ }
\end{enumerate}
We impose the following ``weak isospectral assumption'': 
\begin{enumerate}
\item[(H$_2$)]
{\em There is  $a>0$ such that  
$\forall\,  t \in [0,1]\, ,\  {\rm Spec}\left(\Delta_t
\right)\, \cap [a,+\infty)\  
\subset \   {\mathcal I} \, ,$} where ${\mathcal I}$ is given by (H$_1$). 
\end{enumerate}
Using the asymptotics of the eigenvalues $\lambda_j$ as $j\to \infty$  
we have shown  in \cite{PT3} that the condition (H$_1$)-(H$_2$) is 
``natural''  for any  $d> n/2$ ($n={\rm dim}\, X$), which means  that the usual isospectral 
assumption (\ref{e:isospectral})   implies (H$_1$)-(H$_2$) for any such $d$ and any $c>0$. 
\begin{Th}\hspace{-2mm}{\bf .}\label{Th:main-isospectral}
Let $(X,g)$ be a $3$-dimensional analytic Liouville billiard table of classical type such that the boundary
$\Gamma$ has at least one non-periodic geodesic. Let 
\[
[0,1]\ni t \mapsto K_t\in C^\infty(\Gamma,\R) 
\]
 be a continuous family of real-valued functions 
on $\Gamma$ satisfying the isospectral condition (H$_1$)-(H$_2$).  Suppose that $K_0$ and $K_1$ are
invariant with respect to the group of symmetries $G=(\Z/2\Z)^3$ of $\Gamma$. Then $K_0\equiv K_1$. 
\end{Th}
A similar result has been proved in \cite{PT3} for smooth 2-dimensional billiard tables. The idea of the proof of Theorem \ref{Th:main-isospectral} is as follows. Fix the continuous function $\mu$ by $\mu(\xi)=\langle\pi^+(\xi),n_g\rangle^{-1}$.
First, using  \cite[Theorem 1.1]{PT3} we obtain that 
\begin{equation}
\label{robin-invariant}
{\mathcal R}_{K_1,\mu}(\Lambda)={\mathcal R}_{K_0,\mu}(\Lambda)\,,
\end{equation}
for any Liouville torus $\Lambda$ of a frequency vector satisfying a suitable Diophantine condition. 
Next, we prove that the union of such tori is dense in the union of the two charts $U_j$, $j=1,2$, of ``action-angle'' coordinates in  $B^\ast\Gamma$, which implies (\ref{robin-invariant}) for any torus 
$\Lambda\in {\mathcal F}_b$. Now the claim follows from Theorem \ref{Th:R-rigid}.  
In the same way we prove 
Theorem \ref{Th:Problem_A}. First we obtain that ${\mathcal R}_{K,\mu}(\Lambda)=0$ for a set of ``rational tori'' $\Lambda$. Then we prove that the union of these tori  is dense in $U_1\cup U_2$, and  we apply  Theorem \ref{Th:R-rigid}. We point out that the proof of 
Theorem \ref{Th:main-isospectral} presented in Sect.  \ref{sec:proofs} requires only finite smoothness of $K_t$ (see Theorem \ref{Th:isospectral}). 

An important ingredient in the proof of both theorems is  the density of the corresponding families of invariant tori
in $U_j$, $j=1,2$.  
This follows from  the non-degeneracy of the frequency map for Liouville billiard tables of classical type
studied in Sect.\,\ref{sec:frequencies}. Recall that in any chart $U_j$ of  action-angles coordinates the frequency map
assigns to any value of the momentum map the frequency vector of the minimal power $B^m: U_j\to U_j$, $m\ge 1$,   that leaves invariant
the corresponding Liouville tori $\Lambda\subset U_j$. The frequency map is said to be non-degenerate in $U_j$ if
its Hessian with respect to the action variables is non-degenerate in a dense subset of $U_j$. 
We are interested in the following problem:

\vspace{0.3cm}

\noindent{\large\bf Problem C.} Is the frequency map  non-degenerate in  any chart of action-angle coordinates?

\vspace{0.3cm}

We prove in Theorem \ref{Th:frequencies} that this is true in the charts $U_j$, $j=1,2$, for any  {\em analytic L.B.T. of classical type for which the boundary $\Gamma$ admits at least one non-closed geodesic}. The $3$-axial ellipsoid and more generally any
billiard table  of the two-parameter family of L.B.T.s of classical type  of constant scalar curvature described in 
\cite[Theorem 3]{PT2} has these properties.

The non-degeneracy of the frequency map appears also as a hypothesis  in the Kolmogorov-Arnold-Moser theorem. 
In particular,  Theorem \ref{Th:frequencies} allows us to apply the KAM theorem for the billiard ball maps
associated with small perturbations of the L.B.T.s in \cite[Theorem 3]{PT2}. 
It is a difficult problem to prove that the frequency map of a specific completely integrable system is non-degenerate. 
The non-degeneracy of the frequency map  of completely integrable Hamiltonian systems has been systematically
investigated in \cite{Knorrer}. The main idea in  \cite{Knorrer} is to investigate the system at the singularities of the momentum map. In our case we reduce the system  at the boundary $S^\ast\Gamma$ of $B^\ast\Gamma$. 
To our best knowledge this problem has not been rigorously studied for completely integrable billiard tables even
in the case of the billiard table  associated with the interior of the ellipsoid. 

The article is organized as follows. In Sect. \ref{sec:set-up} we recall certain facts about the billiard ball map and
define a Radon transform for completely integrable billiard tables. Sect. \ref{sec:LBT} is concerned with the construction
of L.B.T.s. First we consider  a cylinder $C = \T_{\om_1}\times\T_{\om_2}\times[-N,N]$, where $\T_{l} = \R/l\Z$ for $l>0$
and $N>0$ and define a ``metric'' $g$ and two Poisson commuting quadratic with respect to the impulses integrals $I_1$ and
$I_2$ of $g$ in $C$.
The non-negative quadratic form $g$ is {\em degenerate} at a submanifold $S$ of $C$. 
To make $g$ a Riemannian metric we consider its push-forward on the  quotient  
$\sigma:C\to\tilde C$ of $C$ with respect to the group generated by two commuting involutions $\sigma_1$ and $\sigma_2$ whose fix point set is just $S$. The main result in this section is Proposition  \ref{Prop:Factorization} which provides $\tilde C$ with a differentiable structure such that the push-forwards  $\tilde g:=\sigma_\ast g$, $\tilde I_1:=\sigma_\ast I_1$ and  
$\tilde I_2:=\sigma_\ast I_2$ are smooth forms, $\tilde g$ is a Riemannian metric on $\tilde C$ and $\tilde I_1$ and 
$\tilde I_2$ are Poisson commuting integrals of $\tilde g$.  
In Sect.  \ref{subsec:parameter} we write an explicit parameterization of the regular tori by means of the values of the momentum map corresponding to the integrals $\tilde I_1$ and $\tilde I_2$. The injectivity of the Radon transform is investigated in Sect. \ref{sec:R-rigidity}. The non-degeneracy of the frequency map of an analytic L.B.T. is investigated in Sect. \ref{sec:frequencies}. 
The proof of Theorem \ref{Th:Problem_A} and Theorem \ref{Th:main-isospectral} is given in Sect. \ref{sec:proofs}. 
In the Appendix we investigate the frequency map and the action-angle coordinates of
completely integrable billiard tables and derive a formula for the frequency vectors of $B^m$.

\section{Invariant manifolds, Leray form, and Radon transform}\label{sec:set-up}
\setcounter{equation}{0} 

In the present section we define the Radon transform for
integrable billiard tables.
First we recall the definition of the billiard ball map $B$ associated to a billiard table
$(X,g)$, ${\rm dim}\, X=n$, with boundary $\Gamma$. Denote by $H\in C^\infty(T^*X,\R)$ the Hamiltonian corresponding to
 the Riemannian metric $g$ on $X$ via the Legendre transformation and set  
\[
S^\ast X : =  \{\xi\in T^\ast X\,:\, H(\xi) =  1\}\, ,\quad 
S^\ast X|_{\Gamma}:= \{\xi\in S^\ast X\,:\, \pi_X(\xi)\in
\Gamma\}\, ,
\]
\[
{S^\ast_\pm}X|_\Gamma:= \{\xi\in S^\ast X|_{\Gamma}\,:\,
\pm \langle \xi,n_g\rangle >0\}\, ,
\]
$n_g$ being the inward unit normal to $\Gamma$.
Denote by $r : T^\ast X|_{\Gamma}\to T^\ast X|_{\Gamma}$ the ``reflection'' at the boundary
given by $r : v\mapsto w$,  where $w|_{T_y\Gamma}=v|_{T_y\Gamma}$ and $\langle w,n_g\rangle + \langle v,n_g\rangle=0$. 
Obviously $r : S^\ast X|_{\Gamma}\to S^\ast X|_{\Gamma}$. 
Take $u\in{S^\ast_+} X|_\Gamma\subset T^*X$ and consider the integral curve $\gamma(t;u)$ of
the Hamiltonian vector field $X_H$ on $T^*X$ starting at $u$. 
If it intersects transversally $S^\ast X|_{\Gamma}$ at a time $t_1>0$ and lies entirely in the   
interior of $S^\ast  X$ for $t\in (0,t_1)$, we set ${\mathcal B}_0(u):=\gamma(t_1,u)\in{S^\ast_-}X|_\Gamma$. 
The set ${\cal O}\subseteq{S^\ast_+} X|_\Gamma$ of all such $u$ is  open in  ${S^\ast_+}X|_\Gamma$. 
The billiard ball map is defined by
\[
{\mathcal B}:=r\circ{\mathcal B}_0: {\cal O}\to{S^\ast_+}X|_\Gamma\,. 
\]
Denote by $B^\ast \Gamma:= \{\xi\in T^\ast\Gamma\,:\, H(\xi) < 1\}$ the (open) coball bundle of $\Gamma$.
The natural projection $\pi_+ : {S^\ast_+}X|_\Gamma\rightarrow B^\ast \Gamma$ assigning to each
$u\in S^\ast X|_{\Gamma}$ the covector $u|_{T_x\Gamma}\in B^*\Gamma$ admits a smooth inverse map
$\pi^{+}: B ^\ast \Gamma \rightarrow{S^\ast_+}X|_\Gamma$. 
The map $B:=\pi_+\circ{\mathcal B}\circ\pi^+$ is defined in the open subset $\pi_+(\cal O)$ of the coball bundle of $\Gamma$ and
it is a smooth symplectic map, i.e. it preserves the canonical symplectic two-form 
$\om=dp\wedge dx$ on $B^*\Ga$. The map $B$ will be called a {\em billiard ball map} as well. 

From now on we assume that the billiard ball map $B : B^*\Ga\to B^*\Ga$ is  globally defined and
{\em completely integrable}. 
By definition\footnote{This is one of the many definitions of complete integrability
of billiard ball map.}, the complete integrability of the billiard ball map of $(X,g)$
means that there exist $n-1$ invariant with respect to $B$ smooth functions
$F_1,...,F_{n-1}$ on $B^*\Ga$ which are functionally independent and in involution with respect
to the canonical Poisson bracket on $T^*\Gamma$, i.e. 
\[
\{F_i,F_j\}=0,\;\;\; 1\le i,j\le n-1.
\]
The functions $F_1,...,F_{n-1}$ are said to be functionally independent in $B^*\Ga$ if the form
$dF_1\wedge...\wedge dF_{n-1}$ does not vanish almost everywhere. 
A function $f$ on $B^*\Ga$ is said to be invariant with respect to the billiard ball map $B$ if $B^*f=f$. 
The invariant functions with respect to the billiard ball map are called also {\em integrals}.
In particular, as $F_1,...,F_{n-1}$ are integrals, then  any non-empty level set
\[
L_c\eqdef\{\xi\in B^*\Ga\, :\, F_1(\xi)=c_1,...,F_{n-1}(\xi)=c_{n-1}\}\, ,\ c=(c_1,\ldots,c_{n-1})\in \R^{n-1}, 
\]
is  invariant with respect to the billiard ball map $B : B^*\Ga\to B^*\Ga$.
By Arnold-Liouville theorem any regular compact component $\Lambda_c$ of $L_c$ is
diffeomorphic to the $(n-1)$-dimensional torus $\T^{n-1}$ and there exists a tubular neighborhood of $\Lambda_c$
in $B^*\Ga$ symplectically diffeomorphic to $\D^{n-1}_r\times\T^{n-1}$ that is supplied with the canonical symplectic
structure $\sum_{k=1}^{n-1} dJ_k\wedge d\theta_k$. Here $\D^{n-1}_r:=\{J=(J_1,...,J_{n-1})\in\R^{n-1}\,:\,|J|<r\}$
for some $r>0$, $\theta=(\theta_1,...,\theta_{n-1})$ are the periodic coordinates on $\T^{n-1}$,
and $|\cdot|$ is the Euclidean norm in $\R^{n-1}$. The coordinates $(J,\theta)$ are called {\em action-angle} coordinates
of the billiard ball map. Recall that $\Lambda_c$ is {\em regular} if  the $(n-1)$-form $dF_1\wedge...\wedge dF_{n-1}$
does not vanish at the points of $\Lambda_c$. Any regular torus $\Lambda_c$ is a Lagrangian submanifold of
$B^*\Ga$ and it is  also called a {\em Liouville torus}. 
 
Assume that the Liouville torus $\Lambda_c$ is invariant with respect to $B^m$ for some  $m\ge 1$, i.e.
$B^m(\Lambda_c) = \Lambda_c$. Let $\al_c$ be a $(n-1)$-form defined in a tubular neighborhood of $\Lambda_c$ in
$B^*\Ga$ so that
\begin{equation}\label{LerayForm}
\om^{n-1}:= \om\wedge...\wedge\om=\al_c\wedge dF_1\wedge...\wedge dF_{n-1}\, .
\end{equation}
It follows from (\ref{LerayForm}) that the restriction $\la_c\eqdef\al_c|_{\Lambda_c}$ of
$\al_c$ to $\Lambda_c$ is uniquely defined. The form $\la_c$ is
a volume form on $\Lambda_c$  which is called {\em Leray form}.
As  $B$ preserves both the symplectic structure $\om$ and the functions $F_1,...,F_{n-1}$,
one obtains from (\ref{LerayForm}) that the restriction of $B^m$ to $\Lambda_c$ preserves $\la_c$. 

Fix a positive continuous function $\mu$ on $B^\ast\Gamma$ and denote by ${\mathcal F}$ the set of all Liouville tori. 
For any continuous function $K$ on $\Gamma$ the mapping ${\mathcal R}_{K,\mu} : {\mathcal F}\to \R$, given by
\begin{equation}\label{Radon}
{\mathcal R}_{K,\mu}\left(\Lambda_c\right) \, :=\, \frac{1}{\la_c(\Lambda_c)}\, \int_{\Lambda_c}(\pi_\Ga^*K)\mu\,\la_c\,,
\end{equation}
is called a {\em Radon transform} of $K$.
It is easy to see that the Radon transform does not depend on the different choices made in the definition of
the Leray form. 
\begin{Remark}\hspace{-2mm}{\bf .}\label{rem:def2}
An alternative definition of the Radon transform would be
\begin{equation}\label{e:def2}
\widetilde{\mathcal R}_{K,\mu}\left(\Lambda_c\right) \, :=\, \frac{1}{\la_c(\Lambda_c)}\, 
\sum_{j=0}^{m-1}\int_{\Lambda_c}(B^*)^j\Big((\pi_\Ga^*K)\mu\Big)\,\la_c
\end{equation}
where $m\ge 1$ is the minimal power of $B$ that leaves $\Lambda_c$ invariant, i.e., $B^m(\Lambda_c)=\Lambda_c$.
Note that \eqref{e:def2} appears as a spectral invariant of \eqref{thespectrum} in \cite{PT3}.
We show in Sect. \ref{sec:frequencies} that for L.B.T. of classical type $m=1$ in the charts $U_1$ and $U_2$.
In particular, \eqref{Radon} and \eqref{e:def2} coincide in this case.
\end{Remark}

There is another notion of complete integrability which is related to the ``billiard flow'' of the billiard
table $(X,g)$ (cf. Definition \ref{Def:integrability_billiard_flow}). 
We reformulate Definition \ref{Def:integrability_billiard_flow} in terms of the cotangent bundle $T^*X$:
A billiard table is {\em completely integrable} if there exist $n$ smooth functions
$H_1,...,H_{n-1}, H_n=H$ in a neighborhood $U$ of  $S^*X$ in $T^*X$ with the following properties:
\begin{itemize}
\item[(i)]
the functions $H_j$ are in involution in $U$ with respect
to the canonical Poisson bracket on $T^*X$, i.e. 
$\{H_i,H_j\}=0,\;\;\; 1\le i,j\le n, $
\item[(ii)]
$H_1,..., H_n$ are functionally independent in  $U$,  
\item[(iii)]
$r^\ast H_j=H_j$ in $U|_\Gamma$ for   $1\le j\le n$. 
\end{itemize}
The properties (i) and (iii) imply that $H_j$ is invariant with respect to the billiard flow in $U$ for any $1\le j\le n$. 
In particular, the functions $F_j=H_j\circ\pi^+$, $1\le j\le n-1$,  are integrals of the billiard ball map $B$.
As $H_1,..., H_n$ are functionally independent in $U$ the billiard ball map is completely integrable if, for example,
the integrals $H_j$ are homogeneous functions with respect to the standard action of $\R^\ast:=\R\setminus 0$ on 
the fibers of $T^\ast X\setminus0$.  In this way we see that the billiard ball map of a completely integrable billiard
table is completely integrable if the integrals are homogeneous functions on the fibers of $T^\ast X\setminus 0$.  

\begin{Def}\hspace{-2mm}{\bf .}\label{Def:R-rigidity}
A billiard table $(X,g)$ with a completely integrable billiard ball map will be called {\em ${\mathcal R}$-rigid} with
respect to the density $\mu$ if Problem B has a positive solution.
\end{Def}

\section{Liouville billiard tables}\label{sec:LBT}
\setcounter{equation}{0} 

\subsection{Construction of Liouville billiard tables}\label{sec:construction}
In this section we describe a class of 3-dimensional completely integrable billiard tables called
{\em Liouville billiard tables}. The interior of an ellipsoid is a particular case of a Liouville billiard table --
see \S\,\ref{subsec:ellipsoid} below as well as \S\,5.3 in \cite{PT2} for the general construction of
Liouville billiard tables of arbitrary dimension, where the integrability of the billiard ball map was deduced
from geodesically equivalence principle. Here we write explicitly integrals of the billiard flow  of a Liouville billiard table which 
are quadratic forms in momenta, and hence, homogeneous functions of degree $2$ on the fibers of $T^\ast X\setminus0$.  

For any $N>0$ and any $\om_k>0$ $(k=1,2)$ consider the cylinder 
\[ 
C\eqdef\{(\te_1,\te_2,\te_3)\}\cong\T_{\om_1}\times\T_{\om_2}\times[-N,N]\,,\,\,\,\,\,\,\T_l:=\R/l\;\Z\,,
\] 
where $\te_1$ and $\te_2$ are periodic coordinates with minimal periods $\om_1$ and $\om_2$ respectively and
$\te_{3}$ takes its values in the closed interval $[-N,N]$. 
Define the involutions $\si_1,\si_{2} : C\to C$ of the cylinder $C$ by 
\begin{equation}\label{e:si_1} 
\si_1 : (\te_1,\te_2,\te_3)\mapsto 
(-\te_1,\frac{\om_2}{2}-\te_2,\te_3) 
\end{equation} 
and
\begin{equation}\label{e:si_2} 
\si_2 : (\te_1,\te_2,\te_3)\mapsto 
(\te_1,-\te_2,-\te_3)\,. 
\end{equation} 
As the commutator $[\si_1,\si_2]$ vanishes one can define the action of the Abelian group
${\cal A}:=\Z_2\oplus\Z_2$ on $C$ by  $(\al,\te)\mapsto \alpha\cdot\theta :=(\si_1^{\al_1}\circ\si_2^{\al_2})(\te),$
where $\al=(\al_1,\al_2)\in{\cal A}$ and $\te\in C$. 
Consider the equivalence relation $\sim_{\cal A}$ on $C$ defined as follows: 
The points $p,q\in C$ are equivalent $p\sim_{\cal A} q$ iff they belong to 
the same orbit of ${\cal A}$ (i.e., there is  $\al\in{\cal A}$ such that 
$\al\cdot p=q$). Denote by ${\tilde C}$ the topological quotient $C/\sim_{\cal A}$ of $C$
with respect to the action of ${\cal A}$ and let 
\begin{equation}\label{e:si}
\si : C\to{\tilde C}
\end{equation}
be the corresponding projection.
A point $p\in C$ is called a {\em regular} point of the projection \eqref{e:si}
iff it is not a fixed point of the action for any $0\ne\al\in{\cal A}$. 
The points in $C$ that are not regular will be  called {\em singular} or 
{\em branched points} of the projection $\si$. 
The set of singular points is given by $S\eqdef S_1\sqcup S_2$,  where
\begin{equation*} 
S_1\eqdef
\Big\{\Big(\te_1\equiv 0\;\Big(\mod\frac{\om_1}{2}\Big),\te_2\equiv\frac{\om_2}{4}\;\Big(\mod\frac{\om_2}{2}\Big),\te_3\Big)
\;:\;-N\le \te_3\le N\Big\}
\end{equation*}
and 
\begin{equation*}
S_2\eqdef\Big\{\Big(\te_1,\te_2\equiv 0\;\Big(\mod\frac{\om_2}{2}\Big),\te_3=0\Big)\;:\;\te_1\in\T_{\om_1}\Big\}\,. 
\end{equation*} 
The set $S_1\subset C$ has four connected components homeomorphic to the unit interval $[0,1]\subset\R$
while  $S_2\subset C$ has two connected components homeomorphic to $\T$.
\begin{Lemma}\hspace{-2mm}{\bf .} \label{lem:factor}
The space $\tilde C$ is homeomorphic to the unit disk ${\D}^3$ in $\R^3$. 
The map $\si : C\to{\tilde C}$ is a $4$-folded branched covering of 
$\tilde C$. 
\end{Lemma}
\begin{Remark}\hspace{-2mm}{\bf .}
The image ${\tilde S}_1$ of $S_1$ under the projection $\si : C\to{\tilde C}$ is homeomorphic to the disjoint
union of two unit intervals and the image  ${\tilde S}_2$ of $S_2$ is homeomorphic to $\T$.
\end{Remark}
{\em Proof of Lemma \ref{lem:factor}.}
First consider the action of the involution $\si_1$ on the cylinder
\[
C=\T_{\om_1}\times\T_{\om_2}\times[-N,N].
\]
For any value $c\in[-N,N]$ the involution $\si_1$ is acting on the $2$-torus 
$\T^2_{c}\eqdef\T_{\om_1}\times\T_{\om_2}\times\{\te_3=c\}$ by
\[
\si_1(c) : (\te_1,\te_2)\mapsto\Big(-\te_1,\frac{\om_2}{2}-\te_2\Big)\,.
\]
The involution $\si_1(c) : \T^2_{c}\to\T^2_{c}$ has four fixed points and it is easy to see
that the topological quotient  $\s^2(c)$ of $\T^2_{c}$ with respect to the orbits of the action of $\si_1(c)$
is homeomorphic to the $2$-sphere $\s^2\eqdef\{x\in\R^3\;: \;|x|^2=1\}$.
Hence,
\begin{equation}\label{e:identification}
C_1\eqdef(C/\sim_{\si_1})\cong\{(\s^2(c),c)\;: \;c\in[-N,N]\}.
\end{equation}
Under the identification \eqref{e:identification}, the involution $\si_2 : C_1\to C_1$ becomes
\[
\begin{array}{lcrr}
\s^2\times[-N,N]&\rightarrow&\s^2\times[-N,N]\\
(x_1,x_2,x_3;c)&\mapsto&(x_1,x_2,-x_3;-c)\, .
\end{array}
\]
The fixed points of this involution form a submanifold, $\{(x_1,x_2,0;0)\;: \;x_1^2+x_2^2=1\}\cong\T$, and
the corresponding quotient is homeomorphic to $\D^3$, hence, ${\tilde C}\cong(C_1/\sim_{\si_2})\cong\D^3$.
\finishproof

In what follows we will define a differential structure ${\cal D}$ on ${\tilde C}$ 
and a smooth Riemannian metrics $\tg$ on the manifold 
$X\eqdef({\tilde C},{\cal D})\cong{\D}^3$ such that the billiard table $(X,\tg)$ becomes completely integrable. 
The branched covering $\si : C\to{\tilde C}$ defined above will play an important 
role in our construction. 
To this end choose three real-valued $C^\infty$-smooth functions
$\varphi_1, \varphi_2 : \R\to\R$ and $\varphi_3 :[-N,N]\to\R$
satisfying the following properties:
\begin{itemize}  
\item[$(A_1)$] $\varphi_k$ $(k=1,2,3)$ is an even function depending only on the variable $\te_k$ and\\
$\varphi_1(\te_1)\ge\varphi_2(\te_2)\ge 0\ge\varphi_3(\te_3)$; \\
$\varphi_k$ $(k=1,2)$ is periodic with period $\om_k$;\\
$\varphi_2$ satisfies the additional 
symmetry $\displaystyle\varphi_2(\te_2)=\varphi_2\Big(\frac{\om_2}{2}-\te_2\Big)$, 
\item[$(A_2)$] 
$\nu_k\eqdef\min\varphi_k=\max\varphi_{k+1}$ $(k=1,2)$
and $\nu_1>\nu_2=0$; \\
for any $k\in\{1,2\}$, $\varphi_k(\te_k)=\nu_k$ iff
$\displaystyle\te_k\equiv 0\,\Big(\mod\frac{\om_k}{2}\Big)$;\\
$\varphi_2(\te_2)=\nu_1$ iff 
$\displaystyle\te_2\equiv\frac{\om_2}{4}\,\Big(\mod\frac{\om_2}{2}\Big)$;\\
$\varphi_3(\te_3)=\nu_2$ iff $\te_3=0$;
\item[$(A_3)$] {\em compatibility conditions}:
\begin{itemize} 
\item[(1)] for $k\in\{1,2\}$, $\displaystyle\varphi_k''(0)=\varphi_k''(\om_k/2)>0$ and
$\varphi_1^{(2l)}(0)=\varphi_1^{(2l)}(\om_1/2)$;\footnote{Item $(A_5)$ $(i)$ in \cite{PT2}
has to be written similarly.}
\item[(2)] for any $l\ge 0$, 
$\displaystyle\varphi_1^{(2l)}(0)=(-1)^l\varphi_2^{(2l)}(\om_2/4)$ 
and $\varphi_2^{(2l)}(0)=(-1)^l\varphi_3^{(2l)}(0)$. 
\end{itemize} 
\end{itemize} 
Consider the following quadratic forms on $TC$ (quadratic on any fiber $T_\theta C$) 
\begin{equation}\label{e:g}
dg^2\eqdef\Pi_1\,d\te_1^2+\Pi_2\,d\te_2^2+\Pi_3\,d\te_3^2
\end{equation} 
and
\begin{equation} 
\begin{array}{rcll} 
dI_1^2&\eqdef&(\varphi_2+\varphi_3)\Pi_1\,d\te_1^2+
(\varphi_1+\varphi_3)\Pi_2\,d\te_2^2+(\varphi_1+\varphi_2)\Pi_3\,d\te_3^2,\\[0.3cm]
dI_2^2&\eqdef&(\varphi_2\varphi_3)\Pi_1\,d\te_1^2+
(\varphi_2\varphi_3)\Pi_2\,d\te_2^2+(\varphi_2\varphi_3)\Pi_3\,d\te_3^2,
\end{array}
\label{e:integrals}                             
\end{equation}
where 
\[
\Pi_1\eqdef(\varphi_1-\varphi_2)(\varphi_1-\varphi_3),\ 
\Pi_2\eqdef(\varphi_1-\varphi_2)(\varphi_2-\varphi_3),\ \mbox{ and}\ 
\Pi_3\eqdef(\varphi_1-\varphi_3)(\varphi_2-\varphi_3).
\]
We say also that the forms above are quadratic forms on $C$.
Notice that $dg^2$ is {\em degenerate}, it vanishes on $S$. 
\begin{Prop}\hspace{-2mm}{\bf .}\label{Prop:Factorization} 
Assume that the functions $\varphi_1$, $\varphi_2$, $\varphi_3$ satisfy 
$(A_1)\div(A_3)$. Then there exists a differential structure ${\cal D}$ on ${\tilde C}$ 
such that the projection $\si : C\to{\tilde C}$ is smooth and $\si$ is a local 
diffeomorphism in the regular points. The push-forwards  $\tg\eqdef\si_*g$ and $\tI_k\eqdef\si_*I_k$ $(k=1,2)$
are smooth quadratic forms and $\tg$ is a Riemannian metric on $X:=({\tilde C},{\cal D})$.
In addition, the billiard table $(X,\tg)$ is completely integrable and the quadratic forms $\tI_1$, $\tI_2$, and
$\tI_3(\xi):=\tg(\xi,\xi)/2$, considered as functions on $TX$ are functionally independent and
Poisson commuting\footnote{The canonical symplectic structure on $T^*X$
induces a symplectic structure on $TX$ by identifying  vectors and covectors by means of
the Riemannian metric $\tg$.} integrals of the billiard flow of $\tg$. 
\end{Prop} 
{\em Proof of Proposition \ref{Prop:Factorization}.}
Consider the set $S=S_1\sqcup S_2\subset C$ of branched points of the covering 
$\si : C\to{\tilde C}$. Take a point $p=(\te_1^0,\te_2^0,\te_3^0)\in S_1$ and assume
for example that $\te_1^0 = 0$, $\te_2^0 =\frac{\om_2}{4}$ and $\te_3^0\in [-N,N]$. 
Define a new chart $V_1=\{(x_1,x_2,x_3)\}$ in a neighborhood of $p$ by
$x_k:=\te_k-\te_k^0$, $k=1,2$, and $x_3:=\te_3$, where $|x_k|<\om_k/8$ for $k=1,2$ and $|x_3|\le N$. 
In this chart  $p=(0,0,\te_3^0)$ and
\begin{equation}\label{e:g_centered}
dg^2=Q_1\,dx_1^2+Q_2\,dx_2^2+Q_3\,dx_3^2,
\end{equation} 
\begin{equation}\label{e:I1_centered}
dI_1^2=(\phi_2+\phi_3)Q_1\,dx_1^2+
(\phi_1+\phi_3)Q_2\,dx_2^2+(\phi_1+\phi_2)Q_3\,dx_3^2,
\end{equation}
\begin{equation}\label{e:I2_centered}
dI_1^2=(\phi_2\phi_3)Q_1\,dx_1^2+
(\phi_1\phi_3)Q_2\,dx_2^2+(\phi_1\phi_2)Q_3\,dx_3^2,
\end{equation}
where $Q_1\eqdef(\phi_1-\phi_2)(\phi_1-\phi_3)$,
$Q_2\eqdef(\phi_1-\phi_2)(\phi_2-\phi_3)$,
$Q_3\eqdef(\phi_1-\phi_3)(\phi_2-\phi_3)$, $\phi_k(x_k)=\varphi_k(\te_k^0+x_k)$, $k=1,2$,
and $\phi_3(x_3)=\varphi_3(x_3)$. 
Note that $V_1$ is a tubular neighborhood of the chosen component of $S_1$ and it
does not intersect the other components of $S$.
It follows from $(A_1)\div(A_3)$ that the functions $\phi_1$, $\phi_2$, 
and $\phi_3$ are smooth and have the following properties in $V_1$: 
\begin{itemize} 
\item[$(L_1)$] $\phi_k$ is even and depends only on the variable $x_k$; 
\item[$(L_2)$] $\phi_1>\phi_3$ and $\phi_2>\phi_3$; 
\item[$(L_3)$] $\phi_1$ and $\phi_2$ satisfy: 
\begin{itemize} 
\item[$(i)$] $\phi_1>\nu_1$ if $x_1\ne 0$ and  $\phi_1(0)=\nu_1$, $\phi_1''(0)>0$;
\item[$(ii)$] $\phi_2<\nu_1$ if $x_2\ne 0$ and  $\phi_2(0)=\nu_1$, $\phi_2''(0)<0$;
\item[$(iii)$] $\phi_1^{(2l)}(0)=(-1)^l\phi_2^{(2l)}(0)$  for any $ l\ge 0$. 
\end{itemize} 
\end{itemize} 
In the new coordinates, the involution $\si_1|_{V_1}$ becomes $\si_1|_{V_1} : (x_1,x_2,x_3)\mapsto (-x_1,-x_2,x_3)$.
In order to define a differential structure in a neighborhood of $\sigma(p)$ in ${\tilde C}$ consider the mapping
$\Phi_1 : V_1\to {\rm Im}\,\Phi_1 := W_1$,
\begin{equation}\label{e:Phi}
\Phi_1 : (x_1,x_2,x_3)\mapsto(y_1=x_1^2-x_2^2,\;y_2=2x_1x_2,\;y_3=x_3). 
\end{equation}
By Lemma \ref{Lem:smooth} below the push-forwards $\tg|_{W_1}\eqdef{\Phi_1}_*(g|_{V_1})$,
${\tilde I}_1|_{W_1}:={\Phi_1}_*(I_1|_{V_1})$, and  ${\tilde I}_2|_{W_1:}={\Phi_1}_*(I_2|_{V_1})$ are smooth
quadratic forms on $W_1$ and $\tg|_{W_1}$ is positive definite. 
Since $\Phi_1\circ(\si_1|_{V_1})=\Phi_1$ and $\si_2(V_1)\cap V_1=\emptyset$ we can identify
$\sigma|_{V_1}$ with $\Phi_1$ and get a differential structure in the neighborhood of $\si(p)\in{\tilde C}$.
In a similar way we construct a tubular neighborhood $V_2$ of the component
$\te_1^0=\om_2/2$, $\te_2^0 =\frac{\om_2}{4}$ and $\te_3^0\in [-N,N]$ of $S_1$ together with
a mapping $\Phi_2 : V_2\to W_2$ such that the push-forward of $g|_{V_2}$, $I_1|_{V_2}$, and $I_2|_{V_2}$
are smooth quadratic forms on $W_2$.
Consider also the tubular neighborhoods $V_3:=\si_2(V_1)$ and $V_4:=\si_2(V_2)$ of the other two components
of $S_1$ in ${\tilde C}$ together with the mappings
\begin{equation}\label{e:Phi_34}
\Phi_3:=\Phi_1\circ(\si_2|_{V_3}) : V_3\to W_1\;\;\;\;\;\mbox{and}\;\;\;\;\;\;
\Phi_4:=\Phi_2\circ(\si_2|_{V_4}) : V_4\to W_2\,.
\end{equation}
For $j=3,4$ one has $\Phi_j\circ(\si_1|_{V_j})=\Phi_j$, and therefore we can identify $\Phi_j$ with $\si|_{V_j}$.
As the quadratic forms \eqref{e:g} and \eqref{e:integrals} are invariant with
respect to $\si_2$ we obtain from \eqref{e:Phi_34} that ${\Phi_3}_*(g|_{V_3})=\tg|_{W_1}$,
${\Phi_4}_*(g|_{V_4})=\tg|_{W_2}$, ${\Phi_3}_*(I_j|_{V_3})=\tI_j|_{W_1}$, and
${\Phi_4}_*(I_j|_{V_4})=\tI_j|_{W_2}$, $j=1,2$. In particular, the mappings $\si|_{V_3}: V_3\to W_1$ and
$\si|_{V_4}: V_4\to W_2$ and the push-forward of \eqref{e:g} and \eqref{e:integrals} with respect to them
are smooth. Arguing similarly we treat the case $p\in S_2$ and construct a coordinate chart $W_3$ of
${\tilde S}_2=\si(S_2)$ in ${\tilde C}$. 

Covering the image of the branched points of $\si$ by the charts $W_1$, $W_2$, and $W_3$
we get a differential structure on $\sqcup_{j=1}^3 W_j \supset {\tilde S}_1\sqcup{\tilde S}_2$. 
As the set ${\tilde C}\setminus({\tilde S}_1\sqcup{\tilde S}_2)$ consists of regular points of $\si$ we can 
induce a differential structure on it from the differential stricture of the cylinder $C$.
The union of these two differential structures is compatible and defines a differential structure ${\cal D}$
on ${\tilde C}$. Denote by $X$ the smooth manifold $X=({\tilde C},{\cal D})$.
It follows from $(A_1)$ that the forms \eqref{e:g}  and \eqref{e:integrals} on $C$
are invariant under the involutions \eqref{e:si_1} and \eqref{e:si_2}. In particular, 
the push-forwards $\tg:=\si_*g$, ${\tilde I}_1:=\si_*I_1$, and ${\tilde I}_2:=\si_*I_2$ are
smooth quadratic forms on $X\setminus({\tilde S}_1\sqcup{\tilde S}_2)$. Moreover, we have seen that
the push-forwards $\tg$, ${\tilde I}_1$, and ${\tilde I}_2$ are smooth quadratic forms on $W_j$, and 
that $\tg$ is a Riemannian metric in $W_j$ for any $j\in\{1,2,3\}$. 
Hence, the push-forwards $\tg$, ${\tilde I}_1$, and ${\tilde I}_2$ are smooth quadratic forms on $X$ and 
$\tg$ is a Riemannian metric.
We will show that $I_1$ and $I_2$ are integrals of the billiard flow of the metric $g$ on $C\setminus S$.
Indeed, applying the Legendre transformation $p_k=\Pi_k\,{\dot\te}_k$, $k=1,2,3$ (which is well defined 
only on $C\setminus S$) and dropping for simplicity the factor $\frac{1}{2}$ in the Hamiltonian we get
\begin{equation}\label{e:stakel1}
\left\{
\begin{array}{rcll}
H&=&\displaystyle {\frac{1}{\Pi_1}\,p_1^2+\frac{1}{\Pi_2}\,p_2^2+\frac{1}{\Pi_3}\,p_3^2}\\[0.3cm]
I_1&=&\displaystyle{\frac{\varphi_2+\varphi_3}{\Pi_1}\,p_1^2+
\frac{\varphi_1+\varphi_3}{\Pi_2}\,p_2^2+\frac{\varphi_1+\varphi_2}{\Pi_3}\,p_3^2}\\[0.3cm]
I_2&=&\displaystyle{\frac{(\varphi_2\varphi_3)}{\Pi_1}\,p_1^2+\frac{(\varphi_1\varphi_3)}{\Pi_2}\,p_2^2+
\frac{(\varphi_1\varphi_2)}{\Pi_3}\,p_3^2}\\
\end{array}
\right.
\end{equation}
which can be rewritten in {\em St{\"a}kel form} (cf. \cite{Perelomov}, \cite[\S\,2]{TopCrit})
\begin{equation}\label{e:stakel2}
\left\{
\begin{array}{rcll}
p_1^2&=&\varphi_1^2\,H-\varphi_1\,I_1+I_2\\[0.3cm]
p_2^2&=&-\varphi_2^2\,H+\varphi_2\,I_1-I_2\\[0.3cm]
p_3^2&=&\varphi_3^2\,H-\varphi_3\,I_1+I_2\\
\end{array}
\right..
\end{equation}
In particular, the functions $H$, $I_1$, and $I_2$ Poisson commute with respect to the
canonical symplectic form $\om:=dp_1\wedge d\te_1+dp_2\wedge d\te_2+dp_3\wedge d\te_3$ on
the cotangent bundle $T^*(C\setminus S)$ (see for example \cite[Proposition 1]{TopCrit}). 
Moreover, the forms $I_1$ and $I_2$ are invariant
with respect to the reflection map at the boundary $\rho :(TC)|_{\partial C}\to(TC)|_{\partial C}$ given by
\[
(\te_1,\te_2,\pm N,{\dot\te_1},{\dot\te_2},{\dot\te_3})\stackrel{\rho}{\longmapsto}
(\te_1,\te_2,\pm N,{\dot\te_1},{\dot\te_2},-{\dot\te_3})\,.
\]
Hence ${\tilde I}_1$ and ${\tilde I}_2$ are Poisson commuting integrals of the billiard flow of
the metric $\tg$ on $X\setminus\si(S)$. As $\si(S)$ is a $1$-dimensional submanifold in the $3$-manifold
$X$ we get that ${\tilde I}_1$ and ${\tilde I}_2$ are Poisson commuting integrals of the billiard flow of
the metric $\tg$. A direct computation shows that $H$, $I_1$ and $I_2$ in \eqref{e:stakel1} are functionally
independent on $T^*(C\setminus S)$. Hence, ${\tilde H}$, ${\tilde I}_1$ and ${\tilde I}_2$ are
functionally independent on $T^*(X\setminus\si(S))$.
\finishproof 

\begin{Lemma}\hspace{-0.2mm}{\bf .}\label{Lem:smooth}
The quadratic forms $\tg|_{W_1}={\Phi_1}_*(g|_{V_1})$, ${\tilde I}_1|_{W_1}={\Phi_1}_*(I_1|_{V_1})$, and 
${\tilde I}_2|_{W_1}={\Phi_1}_*(I_2|_{V_1})$ are smooth and $\tg|_{W_1}$ is positive definite.
\end{Lemma}
\noindent{\em Proof of Lemma \ref{Lem:smooth}.}
A direct computation involving \eqref{e:Phi} shows that
\[
d\tg^2|_{(W_1\setminus{\tilde S}_1)}=\tg_{11}\,dy_1^2+2\tg_{12}\,dy_1dy_2+\tg_{22}\,dy_2^2+\tg_{33}\,dy_3^2
\]
where
\begin{equation}\label{e:tg_coeff1}
\tg_{11}=\frac{1}{4}\Big(\frac{\phi_1-\phi_2}{x_1^2+x_2^2}\Big)
\left(\frac{(\phi_1-\phi_3)x_1^2+(\phi_2-\phi_3)x_2^2}{x_1^2+x_2^2}\right)\,,
\quad\quad
\tg_{12}=\frac{1}{4}\Big(\frac{\phi_1-\phi_2}{x_1^2+x_2^2}\Big)^2x_1x_2\,,
\end{equation}
and
\begin{equation}\label{e:tg_coeff2}
\tg_{22}=\frac{1}{4}\Big(\frac{\phi_1-\phi_2}{x_1^2+x_2^2}\Big)
\left(\frac{(\phi_1-\phi_3)x_2^2+(\phi_2-\phi_3)x_1^2}{x_1^2+x_2^2}\right)\,,\quad
\tg_{33}=(\phi_1-\phi_3)(\phi_2-\phi_3)\,.
\end{equation}
Let $A:=\phi_1\frac{\p}{\p x_1}\otimes dx_1+\phi_2\frac{\p}{\p x_2}\otimes dx_2+
\phi_1\frac{\p}{\p x_3}\otimes dx_3\in C^\infty(TV_1\otimes T^*V_1)$. 
A similar computation as above shows that 
\[
{\tilde A}|_{(W_1\setminus{\tilde S}_1)}={\tilde A}_{11}\frac{\p}{\p y_1}\otimes dy_1+
{\tilde A}_{12}\frac{\p}{\p y_1}\otimes dy_2+{\tilde A}_{21}\frac{\p}{\p y_2}\otimes dy_1+
\phi_3\frac{\p}{\p y_3}\otimes dy_3
\]
where
\begin{equation}\label{e:tA_coeff}
{\tilde A}_{11}=\frac{\phi_1x_1^2+\phi_2x_2^2}{x_1^2+x_2^2}\,,\,\,
{\tilde A}_{12}=\frac{\phi_1-\phi_2}{x_1^2+x_2^2}\,x_1x_2\,,\,\,
{\tilde A}_{21}=\frac{\phi_1-\phi_2}{x_1^2+x_2^2}\,x_1x_2\,,\,\,
{\tilde A}_{22}=\frac{\phi_1x_2^2+\phi_2x_1^2}{x_1^2+x_2^2}\,.
\end{equation}
Consider the tensor field $A$ as a  section  in ${\rm Hom}\, (T^*V_1,T^*V_1)$. Then we have 
\begin{equation}\label{e:integrals_c}
\det(A+c)g((A+c)^{-1}\xi,\xi)=c^2 g(\xi,\xi)+c I_1(\xi,\xi)+I_2(\xi,\xi) 
\end{equation}
for any  $c>-\max\limits_{0\le\te_1\le\om_1}\varphi_1(\te_1)$. 
We will show that the coefficients \eqref{e:tg_coeff1}, \eqref{e:tg_coeff2}, and \eqref{e:tA_coeff},
when re-expressed in terms of the variables $(y_1,y_2,y_3)$, are smooth in $W_1$.
Then the statement of the Lemma will follow from the relation \eqref{e:integrals_c} 
and the properties of the Vandermonde determinant.

Consider, for example, the function 
\[
\Phi(x_1,x_2):=\frac{\phi_1(x_1)-\phi_2(x_2)}{x_1^2+x_2^2}\ ,\quad (x_1,x_2)\neq (0,0)\, . 
\]
Fix $m\in \N$, $m\ge 1$. 
Using $(L_3)$ and the Taylor formula with an integral reminder term we get 
\begin{eqnarray*}
\phi_1(x_1)&=&\sum_{k=0}^m a_k x_1^{2k} + x_1^{2m+1}S_{1,2m+1}(x_1)\\
\phi_2(x_2)&=&\sum_{k=0}^m (-1)^k a_k x_2^{2k} + x_2^{2m+1}S_{2,2m+1}(x_2)
\end{eqnarray*}
where $S_{j,2m+1}$, $j=1,2$,  are smooth functions in a neighborhood of $0$.  
Lemma \ref{Lem:poynolials} below implies that
\[
\Phi(x_1,x_2)= \sum_{k=0}^{m-1} \Phi_k(y_1,y_2) + S_{2m+1}(x_1,x_2)\quad  \mbox{for}\quad (x_1,x_2)\neq (0,0) \, , 
\]
where $\Phi_k(y_1,y_2):=P_k(y_1,y_2)$ for $k$-odd and  $\Phi_k(y_1,y_2):=R_k(y_1,y_2)$ for $k$-even are homogeneous polynomials of degree $2k$ with respect to $(y_1,y_2)$, and 
\[
S_{2m+1}(x_1,x_2):=\frac{x_1^{2m+1}S_{1,2m+1}(x_1)-x_2^{2m+1}S_{2,2m+1}(x_2)}{x_1^2+x_2^2}\ ,\quad  
(x_1,x_2)\neq (0,0) .
\]
 Consider the directional  derivatives 
\[
\partial_{y_1}:=\frac{\p}{\p y_1}=\frac{1}{2(x_1^2+x_2^2)}\Big(x_1\frac{\p}{\p x_1}-x_2\frac{\p}{\p x_2}\Big)\ \mbox{and}\
\partial_{y_2}:=\frac{\p}{\p y_2}=\frac{1}{2(x_1^2+x_2^2)}\Big(x_2\frac{\p}{\p x_1}+x_1\frac{\p}{\p x_2}\Big).
\]
We have 
\[
\displaystyle \lim_{(x_1,x_2)\to (0,0)}\, \partial_{y_1}^\alpha \partial_{y_2}^\beta\,  S_{2m+1}(x_1,x_2)= 0
\]
for $\alpha+\beta\le m$.    Hence, $\Phi$ can be extended by continuity to a  $C^\infty$-smooth function in the variables $(y_1,y_2)$ in a neighborhood of $(0,0)$ and its Taylor series is $\sum_{k=0}^\infty \Phi_k(y_1,y_2)$. In the case when  $\phi_1$ and $\phi_2$ are real analytic 
the power series $\sum_{k=0}^\infty \Phi_k(y_1,y_2)$ is uniformly convergent in a neighborhood of $(0,0)$. 

Arguing similarly we obtain that the coefficients \eqref{e:tg_coeff1}-\eqref{e:tA_coeff} are $C^\infty$-smooth in the variables $(y_1,y_2)$ when  $\phi_1$ and $\phi_2$ are smooth  and real analytic if $\phi_1$ and $\phi_2$ are real analytic.
Moreover, by Taylor's formula $\phi_1(x_1)=\nu_1+a_1 x_1^2+o(x_1^2)$ as $x_1\to 0$ and
$\phi_2(x_2)=\nu_1-a_1 x_2^2+o(x_2^2)$ as $x_2\to 0$ that together with \eqref{e:tg_coeff1} and \eqref{e:tg_coeff2}
implies  $\tg_{11}=a_1(\nu_1-\phi_3)+o(1)$, $\tg_{12}=o(1)$, and $\tg_{22}=a_1(\nu_1-\phi_3)+o(1)$ as
$y\to((0,0,y_3^0))$. Hence, $d\tg^2|_{(W_1\setminus{\tilde S}_1)}$ can be extended by continuity to
$(0,0,y_3^0)\in{\tilde S}_1$ and by $(L_2)$ the extension is positive definite.
This completes the proof of the Lemma.  
\finishproof

\begin{Lemma}\hspace{-2mm}{\bf .}\label{Lem:poynolials}
For any $m\ge 2$,
\[
x_1^{2m}-x_2^{2m}=
\left\{
\begin{array}{lc}
(x_1^2+x_2^2)\,P_{m-1}(y_1,y_2),&m - \mbox{\rm even},\\
Q_m(y_1,y_2),&m - \mbox{\rm odd},
\end{array}
\right.
\]
\[
x_1^{2m}+x_2^{2m}=
\left\{
\begin{array}{lc}
(x_1^2+x_2^2)\,R_{m-1}(y_1,y_2),&m - \mbox{\rm odd},\\
N_m(y_1,y_2),&m - \mbox{\rm even},
\end{array}
\right.
\]
where $P_m$, $Q_m$, $R_m$, and $N_m$ are polynomials of $y_1$ and $y_2$ of degree $m$.
\end{Lemma}
\noindent{\em Proof of Lemma \ref{Lem:poynolials}.}
Introduce the complex  variables $z:=x_1+i x_2$ and $w:=y_1+i y_2$ and note that $w=z^2$.
Then, for any $m\ge 2$, 
$x_1^{2m}\pm x_2^{2m}=\Big((z+{\bar z})^{2m}\pm(-1)^m(z-{\bar z})^{2m}\Big)/2^{2m}$.
Finally, using Newton's binomial formula one concludes the Lemma.
\finishproof

Following \cite{PT1} we impose the following additional assumptions on the
functions $\varphi_k$:
\begin{itemize}
\item[$(A_4)$]  $\varphi_1(\te_1)=\varphi_1(\om_1/2-\te_1)$
\item[$(A_5)$]  for any $k\in\{1,2\}$ the derivative $\varphi_k'(\te_k)>0$ on $(0,\om_k/4)$ and
$\varphi_3'(\te_3)<0$ on $(0,N]$.
\end{itemize}
The condition  $\varphi_3'(N)<0$ means that the boundary of $X$ is locally geodesically convex.
\begin{Def}\hspace{-2mm}{\bf .}
The billiard table $(X,\tg)$ in Proposition \ref{Prop:Factorization} is
called a {\em Liouville billiard table} (shortly L.B.T.). Liouville billiard tables satisfying
conditions $(A_4)$ and $(A_5)$ are called {\em Liouville billiard tables
of classical type}. In the case when $\varphi_1$, $\varphi_2$, and $\varphi_3$ are real analytic,
the billiard table is called {\em analytic} L.B.T.
\end{Def}
The involutions, 
\begin{eqnarray}
(\te_1,\te_2,\te_3)&\mapsto&(-\te_1,\te_2,\te_3\Big)\nonumber\\
(\te_1,\te_2,\te_3)&\mapsto&\Big(\frac{\om_1}{2}-\te_1,\te_2,\te_3\Big)\label{symmetries}\\
(\te_1,\te_2,\te_3)&\mapsto&(\te_1,-\te_2,\te_3)\nonumber
\end{eqnarray}
induce a group of isometries $G(X)=G(X,\tg)$ on $X$ which is isomorphic to
the direct sum 
\[
G(X)\cong\Z_2\oplus\Z_2\oplus\Z_2\,.
\]
\begin{Remark}\hspace{-2mm}{\bf .}
The action of $G(X)$ on $(X,\tg)$ is an analog of the action of the group
$\Z_2\oplus\Z_2\oplus\Z_2$ in the interior of the ellipsoid in $\R^3=\{(x,y,z)\}$
generated by the reflections with respect to the coordinate planes
$O_{xy}$, $O_{yz}$ and $O_{xz}$.
\end{Remark}
\begin{Remark}\hspace{-2mm}{\bf .} 
The  compatibility conditions 
$\varphi_1^{(2l)}(0)=\varphi_1^{(2l)}(\om_1/2)$, $l=0,1,\ldots $, in $(A_3)$   follows from $(A_4)$ for L.B.T.s of classical type. 
\end{Remark}

\subsection{Ellipsoidal billiard tables}\label{subsec:ellipsoid}
Denote by $\R^3$ the Euclidean space $\R^3=\{(x_1,x_2,x_3)\}$ 
supplied with the standard Euclidean metric  
$dg_0^2\eqdef dx_1^2+dx_2^2+dx_3^2$. 
A class of L.B.T.s in $\R^3$ 
depending on $3$ real parameters $b_1>b_2>b_3$ can be obtained 
using the mapping: 
$$ 
\Sigma_0 : \left\{ 
\begin{array}{l}
 \displaystyle{x_1=\sqrt{(b_1-\la_2)\frac{b_1-\la_3}{b_1-b_3}}\,\cos\phi_1}\\ 
 \displaystyle{x_2=\sqrt{b_2-b_3}\,\sin\phi_1\cos\phi_2}\\ 
 \displaystyle{x_3=\sqrt{\frac{b_3-\la_1}{b_3-b_1}}\,\phi_3\sin\phi_2 }
\end{array} 
\right. 
$$ 
where $\la_k\eqdef b_{k+1}+(b_k-b_{k+1})\sin^2\phi_k$ $(k=1,2)$, 
$\la_3\eqdef b_3-\phi_3^2$, $\phi_k$ $(k=1,2)$ are 
periodic coordinates with period $2\pi$, and $-N\le\phi_3\le N$. 
The mapping $\Sigma_0 : \T^2\times[-N,N]\to\R^3$  gives a $4$-folded branched covering of an ellipsoidal
domain $X$ in $\R^3$ and $(X,g_0)$ is a L.B.T. of classical type -- for details see \S\,5 in \cite{PT2}.
More generally, the two-parameter family of billiard tables $(M^3, g_{\alpha,\beta})$ of constant scalar curvature
$\kappa$ in \cite[Theorem 3]{PT2} consists of L.B.T.s of classical type according to \S\,5.4 in \cite{PT2}.
The boundary of any billiard table of the family is geodesically equivalent to the ellipsoid. In particular,
it has non-periodic geodesics and satisfies the hypothesis of Theorem \ref{Th:Problem_A} and
Theorem \ref{Th:main-isospectral}. This family contains the ellipsoid ($\kappa= 0$) and L.B.T.s of both positive
and negative scalar curvature that are realized on the standard sphere and on the hyperbolic space respectively.

\subsection{Parameterization of the Lagrangian tori}\label{subsec:parameter}
The aim of this section is to obtain charts of action-angle coordinates for L.B.T.s of classical type and to parameterize the corresponding Liouville tori. Recall that 
a L.B.T.  $(X,\tg)$ is obtained as a quotient space
of the cylinder
\[
C=\{(\te_1\,(\mod\om_1),\te_2\,(\mod\om_2),\te_3)\}\cong\T_{\om_1}\times\T_{\om_2}\times[-N,N]
\]
with respect to the group action of ${\cal A}=\Z_2\oplus\Z_2$ as described
in Sect. \ref{sec:construction}. By Proposition \ref{Prop:Factorization},  
the projection $\si : C\to X$ is smooth and invariant with respect to
the group action of ${\cal A}$ on $C$. Moreover, the push-forwards of the quadratic forms (\ref{e:integrals})
with respect to the projection $\si : C\to X$ are integrals of the billiard flow on $(X,\tg)$.
The boundary $\p C$ of $C$ has two connected components defined by $\theta_3=\pm N$ and we set
\[
\T^2_N\eqdef\{(\te_1\,(\mod\om_1),\te_2\,(\mod\om_2),\te_3=N)\}.
\]
By  construction the restriction $\si|_{\T_N^2}$ of the projection $\si : C\to X$ to $\T^2_N$ is
a double branched covering of the boundary $\Ga =\partial X$. 

Denote $C_r:=C\setminus S$ and introduce on $T^*C$ the coordinates $\{(\te_1,\te_2,\te_3;p_1,p_2,p_3)\}$, where
$p_1$, $p_2$, and $p_3$ are the conjugated impulses.
The Legendre transformation corresponding to the Lagrangian $L_g(\xi):=g(\xi,\xi)/2$, $\xi\in T C_r$, 
transforms the Lagrangian and the integrals (\ref{e:integrals}) to the functions $H$, $I_1$ and $I_2$
on $T^*C_r$ given by (\ref{e:stakel1}).\footnote{For simplicity we drop the factor $\frac{1}{2}$ in
the Hamiltonian function.}
Set
\[
Q_1\eqdef  T^*C_r|_{\T^2_N}= \{\eta =(\theta,p)\in T^*C_r\;:\;\te_3=N\}.
\]
The restriction ${\tilde\om}_1$ of the symplectic two-form
\[
\om=dp_1\wedge d\te_1+dp_2\wedge d\te_2+dp_3\wedge d\te_3
\]
to $Q_1$ is ${\tilde\om}_1\eqdef\om|_{Q_1}=dp_1\wedge d\te_1+dp_2\wedge d\te_2$.
This form is degenerate and its kernel $\mathop{\rm Ker}{\tilde\om}_1$ is spanned on the
vector field $\frac{\p}{\p p_3}$.
Denote by $Q$ the {\em isoenergy surface}
\[
Q\eqdef\{\eta\in T^*C_r\;:\;H(\eta)=1\}
\]
and consider the set $Q_2\eqdef Q\cap Q_1$.
It is clear that $Q_2$ is diffeomorphic to the restriction of
the unit cosphere bundle $S^*_gC_r$ of $C_r$ to the torus $\T^2_N$.
The set 
\[
Q_2^+\eqdef\{\eta=(\theta,p)\in Q_2\;:\;p_3<0\}
\]
can be identified with the set $S^*_+C_r|_{\T^2_N}$  of all $\eta$ in $S^*_gC_r|_{\T^2_N}$ such
that $\left<\eta,n_g\right>>0$, where $n_g$ denotes the inward unit normal to
$\T^2_N\setminus S_1$. Moreover, the open coball bundle $B_g^\ast (\T^2_N\setminus S_1)$ can be identified with 
\begin{equation}\label{+condition}
\left\{(\te_1,\te_2;p_1,p_2)\in T^\ast (\T^2\setminus S')\, :\  \left(\frac{p_1^2}{\Pi_1}+\frac{p_2^2}{\Pi_2}\right)|_{\te_3=N} < 1\right\}, 
\end{equation}
where $S' \eqdef\{(\te_1\equiv 0\;(\mod\frac{\om_1}{2}),
\te_2\equiv\frac{\om_2}{4}\;(\mod\frac{\om_2}{2})\}$. 
Consider the  map $R: B^\ast_g (\T^2_N\setminus S_1)\to Q_2^+$  given by
\[
R: (\te_1,\te_2;p_1,p_2) \mapsto (\te_1,\te_2;p_1,p_2,p_3)\, \ \mbox{where}\ p_3=-\sqrt{\Pi_3}\sqrt{1-\frac{p_1^2}{\Pi_1}-\frac{p_2^2}{\Pi_2}}\, .
\]
The coball bundle
$B_g^*(\T^2_N\setminus S_1)$  can be considered as a phase space of the billiard
ball map $B : B^*(\Ga\setminus \sigma(S))\to B^*(\Ga\setminus \sigma(S))$ via the branched double covering $\si|_{\T^2_N} : \T^2_N\to\Ga$. In this setting the map $R$ can be identified with  $\pi^+$. 
We have also 
${\tilde\om}_2\eqdef R^*{\tilde\om}_1=dp_1\wedge d\te_1+dp_2\wedge d\te_2$.
Moreover, the functions $\cI_1\eqdef R^*I_1$ and $\cI_2\eqdef R^*I_2$ are functionally independent integrals of $B$ in $B_g^*(\T^2_N\setminus S_1)$. 

In the coordinates $\{(\te_1,\te_2;p_1,p_2)\}$
the integrals $\cI_1$ and $\cI_2$ become (cf. \eqref{e:stakel1})
\begin{equation}
\cI_1=(\varphi_1+\varphi_2)-(\varphi_1-\nu_3)\frac{p_1^2}{\Pi_1}
-(\varphi_2-\nu_3)\frac{p_2^2}{\Pi_2}\label{cI_1}\, ,
\end{equation}
\begin{equation}
\cI_2=\varphi_1\varphi_2-\varphi_2(\varphi_1-\nu_3)\frac{p_1^2}{\Pi_1}
-\varphi_1(\varphi_2-\nu_3)\frac{p_2^2}{\Pi_2}\label{cI_2}\, ,
\end{equation}
where $\nu_3:=\varphi_3(N)<\nu_2=0$ in view of ($A_1$) and ($A_2$). 
In order to describe the invariant manifolds of the billiard ball map $B$
we choose  real constants $h_1$ and $h_2$ and consider
the level set
\[
{\tL}_h\eqdef\{\cI_1=h_1,\cI_2=h_2\}\subset B_g^*(\T^2_N\setminus S_1),\,\,\,h=(h_1,h_2).
\]
Consider the quadratic polynomial,
\begin{equation}\label{polynomial}
P(t)\eqdef t^2-h_1t+h_2=(t-\ka_1)(t-\ka_2)\, ,
\end{equation}
where $\ka_1$ and $\ka_2$ are the roots of $P$ and $h_1=\ka_1+\ka_2$, $h_2=\ka_1\ka_2$.
If $(\te_1,\te_2;p_1,p_2)\in {\tL}_h$,  it follows from (\ref{e:stakel2}) that 
\begin{equation}\label{p_1}
P(\varphi_1(\te_1))=\varphi_1^2(\te_1)-h_1\varphi_1(\te_1)+h_2= p_1^2\ge 0\, ,
\end{equation}
\begin{equation}\label{p_2}
-P(\varphi_2(\te_2))=-\varphi_2^2(\te_2)+h_1\varphi_2(\te_2)-h_2=p_2^2\ge 0\, ,
\end{equation}
and 
\begin{equation}\label{p_3}
P(\varphi_3(N))=P(\nu_3)= \nu_3^2-h_1\nu_3+h_2\ge 0.
\end{equation}
Then the set ${\tL}_h$ is non-empty if and only if there is a point
$(\te_1,\te_2)\in(\R/\om_1\Z)\times(\R/\om_2\Z)$
such that the inequalities (\ref{p_1}), (\ref{p_2}),
and (\ref{p_3}) are satisfied.
In particular, it follows from (\ref{p_1}) and (\ref{p_2})
that the roots $\ka_1\le\ka_2$
are real, hence,  
${\EuScript D}:=h_1^2-4h_2\ge 0$.  Moreover, $(A_1)\div(A_2)$ imply 
\begin{equation}
\left\{
\begin{array}{lcrr}
\nu_1\le\varphi_1(\te_1)\le\nu_0\, ,\\
0=\nu_2\le\varphi_2(\te_2)\le\nu_1\, ,\\
\nu_3= \varphi_3(N) < \nu_2=0\,.
\end{array}
\right.
               \label{e:varphi}
\end{equation}
Then the  following four cases can occur:
\begin{itemize}
\item[(A)] $\nu_3\le\ka_1\le\nu_2=0$ and $0=\nu_2\le\ka_2\le\nu_1$;
\item[(B)] $\nu_3\le\ka_1\le\nu_2=0$ and $\nu_1\le\ka_2\le\nu_0$;
\item[(C)] $0=\nu_2\le\ka_1\le\ka_2\le\nu_1$;
\item[(D)] $0=\nu_2\le\ka_1\le\nu_1$ and $\nu_1\le\ka_2\le\nu_0$.
\end{itemize}
Consider the union ${\tilde U}_1$ of all ${\tL}_h$ in $B^*(\T^2_N\setminus S_1)$ such that (A)
{\em with strict inequalities} holds for the corresponding  ($\ka_1,\ka_2$).
We will see below that any ${\tL}_h$ in ${\tilde U}_1$ is a disjoint union of Liouville tori. 
In the same way we define $\tilde U_2$ corresponding to (B), $\tilde U_3$ corresponding to (C) and 
$\tilde U_4$ corresponding to (D). Denote $U_j:=\sigma_*(\tilde U_j)\subset B^*\Ga$, $j=1,2,3,4$, where
$\si_*$ is the push-forward of covectors corresponding to $\si : C\to X$. 
\begin{Def}\hspace{-2mm}{\bf .}\label{Def:boundary_cases}
We refer to cases $(A)$ and $(B)$ as to {\em boundary cases} and denote ${\mathcal F}_b:=U_1\cup U_2$.
\end{Def}
\begin{Remark}\hspace{-2mm}{\bf .}\label{Rem:boundary-cases}  
We will see in Sect.\,\ref{sec:frequencies} that the billiard trajectories in $T^\ast X$ issuing from $U_1\cup U_2$ 
``approximate'' the geodesics on the boundary $\Ga$. 
\end{Remark}
We are going to parameterize the invariant tori belonging to the level set  ${\tL}_h$.
To that end we  need the inverse functions of $\varphi_1|_{[0,\om_1/4]}$ and $\varphi_2|_{[0,\om_2/4]}$. 
According to $(A_1)\div(A_5)$  the function  $\varphi_1$ has the following properties. 
It is a  periodic  function of period $\omega_1/2$, 
$\varphi_1(\omega_1/4+\theta_1)=\varphi_1(\omega_1/4-\theta_1)$ for any $\theta_1$,  the map $\varphi_1:[0,\omega_1/4]\longrightarrow [\nu_1,\nu_0]$ is a homeomorphism, 
$\varphi_1'(\theta_1)>0$ in the interval $(0,\omega_1/4)$, and the critical points of $\varphi_1$ at $\theta_1=0$ and $\theta_1=\omega_1/4$ are non degenerate. 
Denote by $f_1:[\nu_1,\nu_0]\longrightarrow [0,\omega_1/4]$ the inverse map of $\varphi_1|_{[0,\om_1/4]}$.
Then $f_1$ is smooth in $(\nu_1,\nu_0)$,  $f_1'>0$ in that interval, and 
\begin{equation}\label{e:f1}
\left\{
\begin{array}{rcll}
f_1(x_1)&=&F^+_1(\sqrt{x_1-\nu_1}) \quad {\rm as} \quad x_1\to\nu_1+0\, ,\\[0.3cm]
f_1(x_1)&=&F^-_1(\sqrt{\nu_0-x_1}) \quad {\rm as} \quad x_1\to\nu_0-0\, ,
\end{array}
\right.
\end{equation}
where $F^\mp_1$ are smooth functions in a neighborhood of $0$, and
\begin{equation}\label{e:F1}
F^+_1(0)=0\, ,\ F^-_1(0)=\omega_1/4\, ,\ (F^+_1)'(0) = \sqrt{2\varphi_1''(0)^{-1}}\ \mbox{and}\ 
(F^-_1)'(0) = - \sqrt{-2\varphi_1''(\omega_1/4)^{-1}}. 
\end{equation}
The function $\varphi_2|_{[0,\om_2/4]}$ has the same properties, and we denote by
$f_2:[0,\nu_1]\longrightarrow [0,\omega_2/4]$ its inverse function. Then $f_2$ is smooth in $(0,\nu_1)$
and $f_2'>0$ in that interval, and 
\begin{equation}\label{e:f2}
\left\{
\begin{array}{rcll}
f_2(x_2)&=&F^+_2(\sqrt{x_2}) \quad {\rm as} \quad x_2\to 0+0\, ,\\[0.3cm]
f_2(x_2)&=&F^-_2(\sqrt{\nu_1-x_2}) \quad {\rm as} \quad x_2\to\nu_1-0\, ,
\end{array}
\right.
\end{equation}
where $F^\mp_2$ are smooth functions in a neighborhood of $0$ and 
\begin{equation}\label{e:F2}
F^+_2(0)=0\, ,\ F^-_2(0)=\omega_2 /4\, ,\ (F^+_2)'(0) = \sqrt{2\varphi_2''(0)^{-1}}\, ,\
(F^-_2)'(0) = - \sqrt{-2\varphi_2''(\omega_2 /4)^{-1}}. 
\end{equation}
Assume that $\tilde L_h\subset\tilde U_1$.
We have $\nu_3<\ka_1<0$ and $0<\ka_2<\nu_1$. It follows from \eqref{p_1}-\eqref{p_2} and \eqref{e:varphi}
that ${\tL}_h$ consists of four connected components $T^{(k)}_h$ $(1\le k\le 4)$ which are diffeomorphic to $\T^2$. 
Moreover,  the image of each  $T^{(k)}_h$  with respect to
the bundle projection $T^*\T^2_N\to\T^2_N$ coincides with
one of the annuli
\[
A^{'}_h\eqdef\{0\le\te_1\le\om_1;\;-f_2(\ka_2)\le\te_2\le f_2(\ka_2)\}
\]
and
\[
A^{''}_h\eqdef\{0\le\te_1\le\om_1;\;\om_2/2-f_2(\ka_2)\le\te_2\le\om_2/2+f_2(\ka_2)\}\, .
\]
Assume that the tori $T^{(1)}_h$ and $T^{(2)}_h$ are projected onto $A^{'}_h$ and similarly,
$T^{(3)}_h$ and $T^{(4)}_h$ are projected onto $A^{''}_h$.
As the map $\si|_{\T^2_N} : \T^2_N\to\Ga$ is invariant with respect to
the involution
\[
\imath : (\te_1,\te_2)\mapsto(-\te_1,\om_2/2-\te_2),
\]
and $\imath(A^{'}_h)=A^{''}_h$, the pairs $(T^{(1)}_h,T^{(2)}_h)$ and
$(T^{(3)}_h,T^{(4)}_h)$ correspond  the same  pair of invariant tori in
$T^*\Ga$, which we identify with $(T^{(1)}_h,T^{(2)}_h)$.
It follows from \eqref{p_1}, \eqref{p_2} and \eqref{e:varphi} that the map
$\displaystyle r_{\ep_1\ep_2}: A^{'}_h \to T^{(1)}_h$ defined by,
\begin{equation}\label{r}
\displaystyle 
(\te_1,\te_2)\stackrel{r_{\ep_1\ep_2}}{\longmapsto}
\left(\te_1,\te_2;\ep_1\sqrt{\varphi_1(\te_1)^2-h_1\varphi_1(\te_1)+h_2},
\ep_2\sqrt{-\varphi_2(\te_2)^2+h_1\varphi_2(\te_2)-h_2}\right)\,,
\end{equation}
gives a parametrization of the torus $T^{(1)}_h$ for $\ep_1=1$ and $\ep_2=\pm 1$.
In the same way, taking $\ep_1=-1$  and $\ep_2=\pm 1$ we parametrize $T^{(2)}_h$.

In the same way one treats the cases (B), (C) and (D). 
In particular, one gets that ${\tilde U}_1$, ${\tilde U}_2$, and ${\tilde U}_3$ have $4$ connected components while
$\tilde U_4$ has $8$ connected components. Similarly, $U_1$, $U_2$, and $U_3$ have $2$ connected components
and $U_4$ has $4$ connected components.

\section{${\mathcal R}$-rigidity}\label{sec:R-rigidity}
\setcounter{equation}{0} 

We are going to prove that 
Liouville billiard tables of classical type are ${\mathcal R}$-rigid with respect to the densities $\mu$ defined
by (\ref{thedensity}). 
\begin{Theorem}\hspace{-2mm}{\bf .}\label{Th:R-rigid}
Let $(X,\tg)$ be a Liouville billiard table of classical type and
let $K\in C(\Ga,\R)$ be invariant with respect to
the action of the group $G(X)$ on $\Ga$. Suppose that ${\mathcal R}_{K,\mu}(\Lambda)=0$ for any
Liouville torus $\Lambda\subset{\mathcal F}_b$. Then $K\equiv 0$. 
\end{Theorem}
\noindent{\em Proof of Theorem \ref{Th:R-rigid}.} First, consider the case when $\mu\equiv 1$.
Denote the pull-back of $K$ under the projection $\sigma|_{\T_N^2}: \T_N^2\to\Gamma$ by ${\mathcal K}$,
${\mathcal K}\in C(\T_N^2)$. Let $\Lambda_h\in{\mathcal F}_b$ be a Liouville torus and let $T_h$ be
a connected component of $(\sigma|_{\T_N^2})^\ast\Lambda_h \subset B^\ast(\T_N^2)$ where
$h=(h_1,h_2)$ are the values of the integrals ${\mathcal I}_1$ and ${\mathcal I}_2$ on $T_h$.
Then we have 
\[
(\la_h(T_h))^{-1}\int_{T_h}{\mathcal K}\la_h= 2 {\mathcal R}_{K,1}(\Lambda_h)=0, 
\]
where $\lambda_h$ is the corresponding Leray's form on $T_h$.
Note that ${\mathcal K}$ is invariant under the involution $(\te_1,\te_2)\mapsto(-\te_1,\frac{\om_2}{2}-\te_2)$ since
$\sigma|_{\T_N^2}$ is invariant under the involution (\ref{e:si_1}) for $\theta_3=N$. 
Recall that that the group $G(X)\cong\Z_2\oplus\Z_2\oplus\Z_2$ defined by  \eqref{symmetries} acts by isometries on $X$
and on its boundary $\Ga$. Since $K$ is invariant under this  action,  the function ${\mathcal K}(\te_1,\te_2)$ is
invariant with respect to the involutions
\begin{equation}\label{isometries1}
(\te_1,\te_2)\mapsto(-\te_1,\te_2),\,\,\,(\te_1,\te_2)\mapsto(\om_1/2-\te_1,\te_2)
\end{equation}
and
\begin{equation}\label{isometries2}
(\te_1,\te_2)\mapsto(\te_1,-\te_2),\,\,\,(\te_1,\te_2)\mapsto(\te_1,\om_2/2-\te_2)\,.
\end{equation}
From now on we consider ${\mathcal K}\in C(\T_N^2,\R)$ which is invariant with respect to the involutions
\eqref{isometries1} and \eqref{isometries2} and such that for any
\begin{equation}\label{radon1}
\int_{T_h}{\mathcal K}\la_h=0\,\quad\forall\,T_h\in{\mathcal F}_b\,.
\end{equation}
First, take $h=(h_1,h_2)$ and assume, for example, that $T_h\subset{\tilde U}_1$.
We shall give an explicit formula for the Leray form on the connected components of $\tilde L_h$,
using the parameterization obtained in Sect. \ref{subsec:parameter}. 
Set 
$T_h:=T^{(1)}_h$, 
and let $T_h^+$ be the ``half torus'' $r_{11}(A^{'}_h)$, where the map $r_{11}$ is defined by (\ref{r}). 
Consider the set
\[
A^{'}_h(\de)\eqdef\{0\le\te_1\le\om_1;\;-f(\ka_2)+\de\le\te_2\le f(\ka_2)-\de\}\, ,
\]
where $\de>0$ is  sufficiently small.
It follows from (\ref{r}) that the functions $(\te_1,\te_2,\cI_1,\cI_2)$
give a coordinate chart in a neighborhood of the branch $T_h^+(\de)\eqdef r_{11}(A^{'}_h(\de))\subseteq T_h^+$.
We will compute the Leray form on it. In the coordinates $\{(\te_1,\te_2;p_1,p_2)\}$ on $B^*\T^2_N$ we have 
\begin{eqnarray}
{\tilde\om_2}\wedge{\tilde\om}_2&=&2\,dp_1\wedge d\te_1\wedge dp_2\wedge d\te_2\nonumber\\[0.5cm]
&=&2\,d(\sqrt{\varphi_1(\te_1)^2-\cI_1\varphi_1(\te_1)+\cI_2})\wedge d\te_1\wedge
d(\sqrt{-\varphi_2(\te_2)^2+\cI_1\varphi_2(\te_2)-\cI_2})\wedge d\te_2\nonumber\\[0.5cm]
&=&-\frac{1}{2}\,\frac{(\varphi_1(\te_1)-\varphi_2(\te_2))\;d\te_1\wedge d\te_2}
{\sqrt{\varphi_1(\te_1)^2-\cI_1\varphi_1(\te_1)+\cI_2}
\sqrt{-\varphi_2(\te_2)^2+\cI_1\varphi_2(\te_2)-\cI_2}}\wedge d\cI_1\wedge d\cI_2\, .\nonumber
\end{eqnarray}
In particular, letting $\de\to 0+0$ we see that the Leray form on $T_h^+$
can be identified with 
\begin{equation}\label{la_h}
\la_h\eqdef\frac{(\varphi_1(\te_1)-\varphi_2(\te_2))\;d\te_1\wedge d\te_2}
{\sqrt{\varphi_1(\te_1)^2-h_1\varphi_1(\te_1)+h_2}
\sqrt{-\varphi_2(\te_2)^2+h_1\varphi_2(\te_2)-h_2}}\, .
\end{equation}
We have
\begin{eqnarray}
\int_{T_h}{\mathcal K}\la_h\nonumber&=&2\int_{-f_2(\ka_2)}^{f_2(\ka_2)}\!\int_0^{\om_1}
\frac{{\mathcal K}(\te_1,\te_2)(\varphi_1(\te_1)-\varphi_2(\te_2))\;d\te_1d\te_2}
{\sqrt{(\varphi_1(\te_1)-\ka_1)(\varphi_1(\te_1)-\ka_2)}
\sqrt{(\varphi_2(\te_2)-\ka_1)(\ka_2-\varphi_2(\te_2))}}\nonumber\\[0.3cm]
&=&16\int_0^{f_2(\ka_2)}\!\!\!\!\int_0^{\om_1/4}
\frac{{\mathcal K}(\te_1,\te_2)(\varphi_1(\te_1)-\varphi_2(\te_2))\;d\te_1d\te_2}
{\sqrt{(\varphi_1(\te_1)-\ka_1)(\varphi_1(\te_1)-\ka_2)}
\sqrt{(\varphi_2(\te_2)-\ka_1)(\ka_2-\varphi_2(\te_2))}}\, ,\nonumber
\end{eqnarray}
as the functions ${\mathcal K}$, $\varphi_1$, and $\varphi_2$ are invariant 
with respect to the involutions (\ref{isometries1}) and (\ref{isometries2}).
Set $\tK(\te_1,\te_2)\eqdef{\mathcal K}(\te_1,\te_2)(\varphi_1(\te_1)-\varphi_2(\te_2))$ and denote 
\begin{equation}\label{M}
M_A(\ka_1,\ka_2)\eqdef\int_0^{f_2(\ka_2)}\!\!\!\!\int_0^{\om_1/4}
\frac{\tK(\te_1,\te_2)\;d\te_1d\te_2}
{\sqrt{(\varphi_1(\te_1)-\ka_1)(\varphi_1(\te_1)-\ka_2)}
\sqrt{(\varphi_2(\te_2)-\ka_1)(\ka_2-\varphi_2(\te_2))}}.
\end{equation}
Then (\ref{radon1}) implies 
\[
M_A(\ka_1,\ka_2) = 0
\]
for any $\ka_1\in(-\nu_3,0)$ and any $\ka_2\in(0,\nu_1)$.
\begin{Remark}\hspace{-2mm}{\bf .}
Note that for any fixed $\ka_2\in(0,\nu_1)$ the function $\ka_1\to M_A(\ka_1,\ka_2)$ can be extended to
an analytic (possibly multivalued) function on $\C\setminus([0,\ka_2]\sqcup[\nu_1,\nu_0])$.
Since it  vanishes  for $\ka_1\in(-\nu_3,0)$ we obtain that
$M_{A}(\ka_2,\ka_1)\equiv 0$, $\forall\ka_1\in\C\setminus([0,\ka_2]\sqcup[\nu_1,\nu_0])$ and
$\forall\ka_2\in(0,\nu_1)$.
\end{Remark}
Set 
\[
\tK_1(x_1,x_2)\eqdef\tK(f_1(x_1),f_2(x_2))f_1'(x_1)f_2'(x_2)\,. 
\]
It follows from (\ref{e:f1}) and  (\ref{e:f2}) that $\tilde K_1\in L^1((\nu_1,\nu_0)\times (0,\nu_1))$. 
More precisely, (\ref{e:f1}) and  (\ref{e:f2}) imply
\begin{Lemma}\hspace{-2mm}{\bf .}\label{L1}
We have 
\begin{equation}\label{e:K}
\tilde K_1(x_1,x_2)=\frac{\tilde K (f_1(x_1),f_2(x_2))F(x_1, x_2)}{\sqrt{x_1-\nu_1}\sqrt{\nu_0-x_1}\sqrt{x_2}\sqrt{\nu_1-x_2}}\, ,
\end{equation}
where the function $(x_1,x_2)\mapsto K (f_1(x_1),f_2(x_2))$ is continuous on $[\nu_1,\nu_0]\times [0,\nu_1]$,  the function 
$F\in C([\nu_1,\nu_0]\times [0,\nu_1])$ does not  dependent on $\tilde K $ and $F>0$. 
\end{Lemma}
Passing to the variables $x_1=\varphi_1(\te_1)$ and
$x_2=\varphi_2(\te_2)$ in (\ref{M}) we get
\begin{equation}\label{McaseA}
M_A(\ka_1,\ka_2)= \int_{\nu_1}^{\nu_0}\!\!\!\!\int_0^{\ka_2}
\frac{\tK_1(x_1,x_2)\;dx_2dx_1}
{\sqrt{(x_1-\ka_1)(x_2-\ka_1)}
\sqrt{(x_1-\ka_2)(\ka_2-x_2)}}\equiv 0
\end{equation}
for any $\ka_1\in(-\infty,0)$ and any $\ka_2\in(0,\nu_1)$. Consider now the case $(B)$. 
Arguing in the same way we obtain 
\begin{equation}\label{McaseB}
M_B(\ka_1,\ka_2):=\frac{1}{16}\int_{T_h}{\mathcal K}\la_h= \int_0^{\nu_1}\!\!\!\!\int_{\ka_2}^{\nu_0}
\frac{\tK_1(x_1,x_2)\;dx_1dx_2}
{\sqrt{(x_1-\ka_1)(x_2-\ka_1)}
\sqrt{(x_1-\ka_2)(\ka_2-x_2)}}\equiv 0\, 
\end{equation}
 for any $\ka_1\in(-\infty,0)$ and  $\ka_2\in(\nu_1,\nu_0)$.
 
In the same way one obtains: \\
\noindent{\sc Case $(C)$:} For any $0<\ka_1 < \ka_2<\nu_1$,
\begin{equation}\label{McaseC}
M_C(\ka_1,\ka_2):=\int_{T_h}{\mathcal K}\la_h=\int_{\ka_1}^{\ka_2}\!\!\!\!\int_{\nu_1}^{\nu_0}
\frac{\tK_1(x_1,x_2)\;dx_1dx_2}
{\sqrt{(x_1-\ka_1)(x_2-\ka_1)}
\sqrt{(x_1-\ka_2)(\ka_2-x_2)}}\,.
\end{equation}

\noindent{\sc Case $(D)$:} For any $\ka_1\in(0,\nu_1)$ and  $\ka_2\in(\nu_1,\nu_0)$,
\begin{equation}\label{McaseD}
M_D(\ka_1,\ka_2):=\int_{T_h}{\mathcal K}\la_h=\int_{\ka_1}^{\nu_1}\!\!\!\!\int_{\ka_2}^{\nu_0}
\frac{\tK_1(x_1,x_2)\;dx_1dx_2}
{\sqrt{(x_1-\ka_1)(x_2-\ka_1)}
\sqrt{(x_1-\ka_2)(\ka_2-x_2)}}\,.
\end{equation}
\begin{Remark}\hspace{-2mm}{\bf .}
In what follows we will not use the identities \eqref{McaseC} and \eqref{McaseD}.
\end{Remark}
Now, we argue as follows: Take a continuous function $\chi$ on the interval $[0,\nu_1]$
and consider the mean
\[
\overline{M}_A(\varphi;\ka_1)\eqdef\int_0^{\nu_1}\!\!\!M_A(\ka_1,\ka_2)\chi(\ka_2)\;d\ka_2\, ,
\]
where $M_A(\ka_1,\ka_2)$ is given by (\ref{McaseA}).
In view of Lemma \ref{L1}, we can apply  Fubini's theorem to the following integral  
\begin{eqnarray}
0&\equiv&\overline{M}_A(\varphi;\ka_1)=\int_0^{\nu_1}\left(\int_0^{\ka_2}\!\!\!\!\int_{\nu_1}^{\nu_0}
\frac{\tK_1(x_1,x_2)\chi(\ka_2)\;dx_1dx_2}
{\sqrt{(x_1-\ka_1)(x_2-\ka_1)}\sqrt{(x_1-\ka_2)(\ka_2-x_2)}}\right)d\ka_2\nonumber\\[0.3cm]
&=&\int_0^{\nu_1}\!\!\!\!\int_{\nu_1}^{\nu_0}\!\!\!\!
\frac{\tK_1(x_1,x_2)}{\sqrt{(x_1-\ka_1)(x_2-\ka_1)}}
\left[\int_{x_2}^{\nu_1}\!\!\!\!\frac{\chi(\ka_2)\;d\ka_2}{\sqrt{(x_1-\ka_2)(\ka_2-x_2)}}\right]
dx_1dx_2\label{A-}
\end{eqnarray}
for any $\ka_1\in(-\infty,0)$ and any $\chi\in C([0,\nu_1])$.  
Similarly, consider the mean
\[
\overline{M}_B(\varphi;\ka_1)\eqdef\int_{\nu_1}^{\nu_0}\!\!\!M_B(\ka_1,\ka_2)\chi(\ka_2)\;d\ka_2
\]
where $\chi$ is a continuous function on the interval $[\nu_1,\nu_0]$. We obtain as above 
\begin{eqnarray}
0&\equiv&\overline{M}_B(\varphi;\ka_1)=\int_{\nu_1}^{\nu_0}\!\!\!\!\int_0^{\nu_1}\!\!\!\!\int_{\ka_2}^{\nu_1}
\frac{\tK_1(x_1,x_2)\chi(\ka_2)\;dx_1dx_2d\ka_2}
{\sqrt{(x_1-\ka_1)(x_2-\ka_1)}\sqrt{(x_1-\ka_2)(\ka_2-x_2)}}\nonumber\\[0.3cm]
&=&\int_0^{\nu_1}\!\!\!\!\int_{\nu_1}^{\nu_0}\!\!\!\!
\frac{\tK_1(x_1,x_2)}{\sqrt{(x_1-\ka_1)(x_2-\ka_1)}}
\left[\int_{\nu_1}^{x_1}\!\!\!\!\frac{\chi(\ka_2)\;d\ka_2}{\sqrt{(x_1-\ka_2)(\ka_2-x_2)}}\right]
dx_1dx_2\label{-B}
\end{eqnarray}
for any $\ka_1\in(-\infty,0)$.
Finally, combining (\ref{A-}) and (\ref{-B})  we obtain  for any $\chi\in C([0,\nu_0])$ and
any $\ka_1\in(-\infty,0)$ the equality 
\begin{equation}
\int_0^{\nu_1}\!\!\!\!\int_{\nu_1}^{\nu_0}\!\!\!\!
\frac{\tK_1(x_1,x_2)}{\sqrt{(x_1-\ka_1)(x_2-\ka_1)}}
\left[\int_{x_2}^{x_1}\!\!\!\!\frac{\chi(\ka_2)\;d\ka_2}{\sqrt{(x_1-\ka_2)(\ka_2-x_2)}}\right]
dx_1dx_2\equiv 0.\label{ABCD}
\end{equation}
In particular, for any $k\ge 0$ and for any $\ka_1\in(-\infty,0)$,
\begin{equation}\label{ABCD_k}
\int_0^{\nu_1}\!\!\!\!\int_{\nu_1}^{\nu_0}\!\!\!\!
\frac{\tK_1(x_1,x_2)}{\sqrt{(x_1-\ka_1)(x_2-\ka_1)}}
R_k(x_1,x_2)\;dx_1dx_2\equiv 0,
\end{equation}
where 
\begin{equation}\label{R_k}
R_k(x_1,x_2)\eqdef\int_{x_2}^{x_1}\!\!\!\!\frac{z^k\;dz}{\sqrt{(x_1-z)(z-x_2)}}=
\int_0^1\frac{(x_2+t(x_1-x_2))^k}{\sqrt{t(1-t)}}\;dt.
\end{equation}
Recall that the Legendre polynomials $P_k$, $k\ge 0$, can be generated by the power series expansion,
\begin{equation}\label{legendre}
(1-2\omega z + z^2)^{-1/2}\, =\, \sum_{k=0}^{\infty}\, P_k(\omega)\, z^k \, ,
\end{equation}
which is convergent for small $z$.  
For $0<x_2\le x_1$ we set $s_1\eqdef(x_1+x_2)/2$ and  $s_2\eqdef\sqrt{x_1x_2}$. 
\begin{Lemma}\hspace{-2mm}{\bf .}\label{R_k=P_k}
For any $k\ge 0$ and for any $0<x_2\le x_1$, $R_k(x_1,x_2) = \pi s_2^k P_k(s_1/s_2)$.
\end{Lemma}
\noindent{\em Proof.} 
For any given values of $x_1$ and $x_2$, $0<x_2\le x_1$, consider the power series in $z$,
\[
I(z)\,:=\,\sum_{k=0}^{\infty}\, R_k(x_1,x_2)\,z^k
\]
There exists $0<r<\infty$ sufficiently small such that the power series converges for $|z|\le r$ and
\[
I(z)\,=\,\int_0^1\,\frac{1}{t}\, \frac{1}{1- z(x_2+t(x_1-x_2))}\,{\sqrt{\frac{t}{1-t}}}\,\,dt\,.
\]
Using the substitution, $s=\sqrt{\frac{t}{1-t}}$ we get
\[
I(z)\,=\, \frac{\pi}{\sqrt{(1-zx_1)(1-zx_2)}}\, ,
\]
and by \eqref{legendre} we obtain
\[
I(z)\, = \, \pi \,  
\sum_{k=0}^{\infty}\,s_2^k P_k(s_1/s_2)\,z^k\,,
\]
which proves the lemma.
\finishproof

Note that the function 
\[
(x_1,x_2)\mapsto Q(x_1,x_2,\kappa_1):=\frac{\tK_1(x_1,x_2)}{\sqrt{(x_1-\ka_1)(x_2-\ka_1)}}
\]
 belongs to $L^1([\nu_1,\nu_0]\times [0,\nu_1])$ in view of Lemma \ref{L1}, and it depends analytically on $\kappa_1\in (-\infty, 0)$. Consider the power series expansion 
\begin{equation}
              \label{e:power-series}
\frac{1}{\sqrt{(x_1-\ka_1)(x_2-\ka_1)}} = \sum_{k=1}^{\infty}Q_j(x_1,x_2)\kappa_1^{-j}\, ,
\end{equation}
where $(x_1,x_2)\in [\nu_1,\nu_0]\times [0,\nu_1]$ and $\kappa_1<0$. Now \eqref{legendre} implies 
$Q_j(x_1,x_2)= -s_2^{j-1} P_{j-1}(s_1/s_2)$
for any $j\ge 1$. 

Using Lemma \ref{R_k=P_k}, (\ref{ABCD_k}) and (\ref{e:power-series}) we obtain that for any $k,j\ge 0$, 
\begin{equation}\label{e:*}
\int_0^{\nu_1}\!\! \int_{\nu_1}^{\nu_0}\,{\tilde K}_1(x_1,x_2)  
s_2^{k+j}P_k(s_1/s_2)P_j(s_1/s_2)\, d x_1 d
x_2 \, =\, 0\, .
\end{equation}
Let $k$ and $m$ be non-negative integers such that $2k\le m$ and let $d$ be the integer part of $m/2$.
We have the following relation due to Adams (see \cite{Adams},
\cite[Chap. XV, Legendre functions, Miscellaneous Examples, Ex. 11]{WW}), 
$$
P_k(z)P_{m-k}(z)\,=\,\sum_{r=0}^{k} c_{k,r}^mP_{m-2r}(z)\,=\,\sum_{r=0}^{d} c_{k,r}^mP_{m-2r}(z)\,,
$$
where for any $0\le r\le d$,
\[
c_{k,r}^m :=\frac{A_{k-r}A_rA_{m-k-r}}{A_{m-r}}\, . \, 
\frac{2m-4r+1}{2m - 2r +1}\,,
\]
with
\[
A_k :=\left\{
\begin{array}{cc}
\displaystyle \frac{1.3.5\ldots (2k-1)}{k},&k\ge 1,\\
1,&k=0,\\
A_{k}=0,&k\le -1\,.
\end{array}
\right.
\]
Hence, for any given $m\ge 0$ we obtain a $(d+1)\times (d+1)$ matrix 
$(c_{k,r}^m)_{k,r=0}^{d}$ which is triangular (all the
elements over the diagonal vanish) and with non-vanishing diagonal elements.
This together with \eqref{e:*} (take $m=j+k$) implies that for any $m\ge 0$, $0\le 2r\le m$, 
$$
\int_0^{\nu_1}\!\!\int_{\nu_1}^{\nu_0}\,{\tilde K}_1(x_1,x_2)  
s_2^{m}P_{m-2r}(s_1/s_2)\, d x_1 dx_2 \, =\, 0\,. 
$$
On the other hand, for any $m\ge 0$ the monomial $z^m$ can be written as a linear combination of
the Legendre polynomials $P_{m-2r}(z)$, $0\le2r\le m$, and we get
\begin{equation}\label{e:A}
\int_0^{\nu_1}\!\! \int_{\nu_1}^{\nu_0}\,{\tilde K}_1(x_1,x_2)  
s_1^{m-2r}s_2^{2r}\, d x_1 dx_2 \, =\, 0\,.
\end{equation}
Consider the set of monomials ${\mathcal M}=\{s_1^{m-2r}s_2^{2r}\,:\,r,m\in\Z,\,0\le 2r\le m\}$.
Obviously $\mathcal M$ is closed under multiplication, $1\in{\mathcal M}$, and it separates the points
$(x_1,x_2)$ of the compact $[\nu_1,\nu_0]\times [0,\nu_1]$, 
since $s_1,s_2^2\in{\mathcal M}$ and  $0\le x_2\le x_1$ are the unique solutions of $x^2-2s_1 x+s_2^2 =0$.
The Stone-Weierstrass theorem implies that the vector space  ${\rm Span}(\mathcal M)$ of all
finite linear combinations of monomials of $\mathcal M$ is dense in $C([\nu_1,\nu_0]\times [0,\nu_1])$.
Choose $\psi\in C([\nu_1,\nu_0]\times [0,\nu_1])$. 
Then for any $\varepsilon>0$ there is $P\in{\rm Span}(\mathcal M)$ such that 
\[
\|P-\psi\|_{C([\nu_1,\nu_0]\times [0,\nu_1])} < \varepsilon\, .
\]
Now (\ref{e:A}) implies 
\[
\begin{array}{lcrr}
\displaystyle\left|\int_0^{\nu_1}\!\!\int_{\nu_1}^{\nu_0}\,{\tilde K}_1(x_1,x_2)\psi(x_1,x_2)d x_1 d x_2\right|\\[0.3cm]
\displaystyle 
\le\int_0^{\nu_1}\!\!\int_{\nu_1}^{\nu_0}\, |{\tilde K}_1(x_1,x_2)|| \psi(x_1,x_2)-P(x_1,x_2)| d x_1 d x_2
\le \varepsilon \|{\tilde K}_1\|_{L^1((\nu_1,\nu_0)\times(0,\nu_1))}
\,.
\end{array}
\]
Hence,
\[
\int_0^{\nu_1}\!\! \int_{\nu_1}^{\nu_0}\,{\tilde K}_1(x_1,x_2)\psi(x_1,x_2) d x_1 d x_2 = 0
\]
for any $\psi \in C([\nu_1,\nu_0]\times [0,\nu_1])$ which implies 
${\tilde K}_1\equiv 0$ on that compact. In particular, ${\mathcal K}\equiv 0$, and hence $K\equiv 0$.
This completes the proof when $\mu\equiv 1$. 

Now, consider the case when $\mu(\xi)=\langle\pi^+(\xi),n_g\rangle^{-1}$.
Assume that $\xi\in\La_h$ where $\La_h$ is a Liouville torus in ${\mathcal F}_b$
and $h=(h_1,h_2)$ are the values of the integrals $\tI_1$ and $\tI_2$ on $\La_h$.
Let $\La_h\subset U_1$. Using the mapping \eqref{r}, we introduce coordinates  $\{(\te_1,\te_2)\}$ on the ``half" tori
$T_{h+}^{(1)}=r_{11}(A_h^{'})$ and $T_{h-}^{(1)}=r_{1,-1}(A_h^{'})$ of $T_h^{(1)}$ 
as well as on $T_{h+}^{(2)}=r_{-11}(A_h^{'})$ and $T_{h-}^{(2)}=r_{-1,-1}(A_h^{'})$ of $T_h^{(2)}$.
Similarly, we parametrize the Liouville tori $\La_h$ in $U_2$.
\begin{Lemma}\hspace{-2mm}{\bf .}\label{Lem:mu}
In coordinates $\{(\te_1,\te_2)\}$ on $\La_h\subset{\mathcal F}_b$, $\mu(\xi)=\langle\pi^+(\xi),n_g\rangle^{-1}$ is
given by
\begin{equation}\label{e:mu}
\mu(\te_1,\te_2)=\frac{\sqrt{(\varphi_1(\te_1)-\nu_3)(\varphi_1(\te_1)-\nu_3)}}
{\sqrt{(\ka_1-\nu_3)(\ka_2-\nu_3)}}\,
\end{equation}
where $\ka_1=h_1+h_1$ and $\ka_2=h_1h_2$.
\end{Lemma}
\noindent{\em Proof of Lemma \ref{Lem:mu}.}
Fix $h=(h_1,h_2)$ so that $\La_h\subset{\mathcal F}_b$.
It follows from \eqref{e:g} that
\[
n_g=-\frac{1}{\sqrt{\Pi_3}}\,\frac{\p}{\p\te_3}\,.
\]
On the other hand, the third equation in \eqref{e:stakel2} shows that
\[
p_3=-\sqrt{\nu_3^2-h_1\nu_3+h_2}=-\sqrt{(\ka_1-\nu_3)(\ka_2-\nu_3)}.
\]
Hence, $\langle\pi^+(\xi),n_g\rangle=\sqrt{(\ka_1-\nu_3)(\ka_2-\nu_3)}/\sqrt{\Pi_3}$.
\finishproof
\begin{Remark}\hspace{-2mm}{\bf .}
The statement of Lemma \ref{Lem:mu} holds also for any $\La_h\in {\mathcal F}$ not necessarily in ${\mathcal F}_b$.
\end{Remark}
Note that the denominator in \eqref{e:mu} is a positive constant on $T_h$ and
the numerator is independent of $h_1$ and $h_2$ and does not vanish. The relation  \eqref{radon1} with
\[
{\tilde K}(\te_1,\te_2)={\mathcal K}(\te_1,\te_1)(\varphi_1(\te_1)-\varphi_2(\te_2))
\sqrt{(\varphi_1(\te_1)-\nu_3)(\varphi_1(\te_1)-\nu_3)} .
\]
 implies that
the expression \eqref{M} 
vanishes. In particular, \eqref{McaseA} and \eqref{McaseB} hold. Finally,
arguing in the same way as in the case $\mu\equiv 1$ one concludes that $K\equiv 0$.
\finishproof

\section{Non-degeneracy of the frequency map}\label{sec:frequencies}
\setcounter{equation}{0} 

In this section we investigate the non-degeneracy of the frequency map of Liouville billiard
tables of classical type. 
\begin{Theorem}\hspace{-2mm}{\bf .}\label{Th:frequencies}
Let $(X,g)$ be an analytic $3$-dimensional Liouville billiard table of classical type. Suppose that there is at least one non-periodic geodesic on $\Gamma$. 
Then the  frequency map  is non-degenerate in the union $U_1\cup U_2$ corresponding to the boundary cases (A) and (B). 
\end{Theorem}
{\em Proof.} As in Sect.\,\ref{subsec:parameter} we introduce coordinates $\{(\te_1,\te_2,\te_3;p_1,p_2,p_3)\}$
on the cotangent bundle $T^*C$,  where $p_1$, $p_2$, $p_3$ are the conjugate variables to
$\te_1$, $\te_2$, and $\te_3$. Solving the system of equations (\ref{e:stakel1}) with respect to
$p_1^2$, $p_2^2$, and $p_3^2$, where $H=1$, $I_1=h_1$, and $I_2=h_2$ are given values of the integrals,  we get 
\begin{equation}\label{impulses}
\left\{
\begin{array}{ccc}
p_1^2&=&\varphi_1^2-h_1\varphi_1+h_2\\
p_2^2&=&-(\varphi_2^2-h_1\varphi_2+h_2)\\
p_3^2&=&\varphi_3^2-h_1\varphi_3+h_2\, .
\end{array}
\right.
\end{equation}
In particular, 
it follows from \eqref{impulses} that the invariant set
\begin{equation}\label{torus}
{\cal T}_{h}:=\{H=1, I_1=h_1, I_2=h_2\}\subset T^*C
\end{equation}
is non-empty if and only if the quadratic polynomial $P(t)=t^2-h_1t+h^2$ has real roots $\ka_1\le\ka_2$ (i.e., ${\cal D}=h_1^2-4h_2\ge 0$).
As in Sect.\,\ref{sec:R-rigidity} we obtain  four cases related to the position of the roots $\ka_1$ and $\ka_2$ with respect to the constants
$\nu_3<\nu_2=0<\nu_1 <\nu_0$, namely,  
\begin{itemize}
\item[(A)] $\nu_3\le\ka_1\le\nu_2=0$ and $0\le\ka_2\le\nu_1$;
\item[(B)] $\nu_3\le\ka_1\le 0$ and $\nu_1\le\ka_2\le\nu_0$;
\item[(C)] $0\le\ka_1\le\ka_2\le\nu_1$;
\item[(D)] $0\le\ka_1\le\nu_1$ and $\nu_1\le\ka_2\le\nu_0$.
\end{itemize}
Recall that $\nu_3=\min\varphi_3$, $\nu_2=\max\varphi_3=\min\varphi_2=0$, 
$\nu_1=\max\varphi_2=\min\varphi_1$, and $\nu_0=\max\varphi_1$. 
In what follows we consider $\ka_1$ and $\ka_2$ as new parameters (constants of motion)\footnote{$h_1=\ka_1+\ka_2$ and $h_2=\ka_1\ka_2$} that
parametrize the invariant set \eqref{torus}.

\vspace{0.5cm}

We first consider the case $(A)$ where $\nu_3\le\ka_1\le\nu_2=0$ and $0=\nu_2\le\ka_2\le\nu_1$.
It follows from \eqref{impulses} that the impulses are real-valued if and only if
\begin{equation}\label{teta}
\left\{
\begin{array}{l}
\nu_1\le\varphi_1(\te_1)\le\nu_0\\
0\le\varphi_2(\te_2)\le\ka_2\\
\nu_3\le\varphi_3(\te_3)\le\ka_1
\end{array}
\right.
\end{equation}
Hence, the projection of the invariant set \eqref{torus} onto the base $C$ is described by
the following inequalities:
\[
0\le\te_1\le\om_1\,;
\]
\[
-f_2(\ka_2)\le\te_2\le f_2(\ka_2)\;\;\;\mbox{or}\;\;\;-f_2(\ka_2)+\frac{\om_2}{2}\le\te_2\le f_2(\ka_2)+\frac{\om_2}{2}\,;
\]
\[
f_3(\ka_1)\le\te_3\le N\;\;\;\mbox{or}\;\;\;-N\le\te_3\le-f_3(\ka_1),
\]
where $f_2$ is the inverse of $\varphi_2|_{[0,\,\om_2/4]}$ and $f_3$ is the inverse of $\varphi_3|_{[0,N]}$.
These inequalities give four rectangular boxes in $C$ that project onto an unique set in ${\tilde C}$
via the projection \eqref{e:si}. Consider, for example, the rectangular box $B_h$ given by
\begin{equation}
B_h : \left\{
\begin{array}{l}
0\le\te_1\le\om_1\,;\\
-f_2(\ka_2)\le\te_2\le f_2(\ka_2)\,;\\
f_3(\ka_1)\le\te_3\le N\,.
\end{array}
\right.
\end{equation}
For any given $\te\in B_h$ we obtain from \eqref{impulses} that
\[
p_1(\te)=\ep_1\sqrt{(\varphi_1(\te_1)-\ka_1)(\varphi_1(\te_1)-\ka_2)}
\]
\[
p_2(\te)=\ep_2\sqrt{(\varphi_2(\te_2)-\ka_1)(\ka_2-\varphi_2(\te_2))}
\]
\[
p_3(\te)=\ep_3\sqrt{(\ka_1-\varphi_3(\te_3))(\ka_2-\varphi_3(\te_3))}
\]
where $\ep_k=\pm 1$. Then the mapping $r_+ : B_h\to T^*C$,
\[
(\te_1,\te_2,\te_3)\mapsto(\te_1,\te_2,\te_3;p_1(\te),p_2(\te),p_3(\te)),
\]
where $\ep_1=1$, $\ep_2=\pm 1$, and $\ep_3=\pm 1$, parametrizes one of the two
connected components of the subset $T_h:=\{{\tilde H}=1, {\tilde I}_1=h_1, {\tilde I}_2=h_2\}\subset T^*X$.
Assume that the strict inequalities $\nu_3<\ka_1<\nu_2=0$ and $0<\ka_2<\nu_1$ hold.

\begin{Remark}\hspace{-2mm}{\bf .}\label{rem:m=1}
This component is diffeomorphic to $\T^2\times[0,1]$ and its intersection with the boundary $T^*X|_\Ga$ of
$T^*X$ has two components which can be identified with the two components of the image of the slice 
\[\{0\le\te_1\le\om_1,\,-f_2(\ka_2)\le\te_2\le f_2(\ka_2),\,\te_3=N\}\]
of $B_h$ with respect to $r_+$ with $\ep_3=1$ and $\ep_3=-1$ respectively. In particular, the impulse
$p_3$ takes constant values of different sign on them. Moreover, the reflection map
$r : T^*X|_\Ga\to T^*X|_\Ga$ is given by
\[
(\te_1,\te_2;p_1,p_2,p_3)\mapsto(\te_1,\te_2;p_1,p_2,-p_3)\,.
\]
Hence, the reflection map interchanges these two components, and by Lemma \ref{Lem:glued_torus} $(c)$, $m=1$
(cf. Remark \ref{rem:def2}). Similarly, we get  $m=1$ in the case $(B)$.
\end{Remark}

Now we compute the generalized actions of the billiard flow corresponding
to $T_h$ (see \eqref{actions}, Appendix),
\begin{eqnarray}
J_1(\ka_1,\ka_2)&=&
\frac{1}{2\pi}\int_0^{\om_1}\!\!\!\sqrt{(\varphi_1(\te_1)-\ka_1)(\varphi_1(\te_1)-\ka_2)}\,\,d\te_1\nonumber\\[0.3cm]
&=&\frac{2}{\pi}\int_{\nu_1}^{\nu_0}\sqrt{(x_1-\ka_1)(x_1-\ka_2)}\,\rho_1(x_1)\,dx_1\,,\label{action1}
\end{eqnarray}
\begin{eqnarray}
J_2(\ka_1,\ka_2)&=&
\frac{1}{\pi}\int_{-f_2(\ka_2)}^{f_2(\ka_2)}\!\!\!\sqrt{(\varphi_2(\te_2)-\ka_1)(\ka_2-\varphi_2(\te_2))}\,\,d\te_2\nonumber\\[0.3cm]
&=&\frac{2}{\pi}\int_0^{\ka_2}\sqrt{(x_2-\ka_1)(\ka_2-x_2)}\,\rho_2(x_2)\,dx_2\,,\label{action2}
\end{eqnarray}
\begin{eqnarray}
J_3(\ka_1,\ka_2)&=&
\frac{1}{\pi}\int_{f_3(\ka_1)}^N\!\!\!\sqrt{(\ka_1-\varphi_3(\te_3))(\ka_2-\varphi_3(\te_3))}\,\,d\te_3\nonumber\\[0.3cm]
&=&\frac{1}{\pi}\int_{\nu_3}^{\ka_1}\sqrt{(\ka_1-x_3)(\ka_2-x_3)}\,\rho_3(x_3)\,dx_3\, ,\label{action3}
\end{eqnarray}
where 
\[
\rho_1(x_1):=f_1'(x_1)>0\, ,\ \rho_2(x_2):=f_2'(x_2)>0\, \ \mbox{and}\ 
\rho_3(x_3):=-f_3'(x_3)>0 
\]
are analytic functions in the intervals $(\nu_1,\nu_0)$,  $(0, \nu_1)$ and $(\nu_3,0)$, respectively, and $f_1$ and $f_2$ satisfy \eqref{e:f1} and  \eqref{e:f2}.  Notice that the functions $F^\pm_1$ and $F^\pm_2$ in \eqref{e:f1} and  \eqref{e:f2} are smooth and even analytic in a neighborhood of $0$, and $f_3(x_3)= \sqrt {-x_3}F_3(\sqrt {-x_3})$, 
where $F_3$ is analytic in a neighborhood of $[\nu_3,0]$. 
By the assumption $(A_5)$, $\rho_3(x_3)$ is smooth 
at $x_3=\nu_3$. In particular, we obtain 
\begin{Remark}\hspace{-2mm}{\bf .}\label{rem:rho} The functions $\rho_1$, $\rho_2$ and $\rho_3$ are analytic in  the intervals $(\nu_1,\nu_0)$,  $(0, \nu_1)$ and $(\nu_3,0)$, respectively, and
\[
\begin{array}{lcrr}
\rho_1(x_1)=\frac{G_1^-(\sqrt{x_1-\nu_1})}{\sqrt{x_1-\nu_1} }\ \mbox{as}\ x_1\to\nu_1+0\ \mbox{and}\ 
\rho_1(x_1)=\frac{G_1^+(\sqrt{\nu_0-x_1})}{\sqrt{\nu_0-x_1}}  \ \mbox{as}\ x_1\to\nu_0 -0\,\\[0.3cm]
\rho_2(x_2)=\frac{G_2^-(\sqrt{x_2})}{\sqrt{x_2}} \ \mbox{as}\ x_2\to 0+0\ \mbox{and}\ 
\rho_2(x_2)=\frac{G_2^+(\sqrt{\nu_1-x_2})}{\sqrt{\nu_1-x_2}}  \ \mbox{as}\ x_2\to\nu_1 -0\,\  \mbox{and}\\[0.3cm]
\rho_3(x_3)=\frac{G_3(\sqrt{-x_3})}{\sqrt{-x_3}} \ \mbox{as}\    x_3\to 0-0,  
\end{array}
\]
where $G_1^\pm$ and $G_2^\pm$ are analytic in a neighborhood of $0$, and $G_3$ is analytic in a neighborhood of $[\nu_3,0]$. Moreover, $G_3(0)>0$ and by \eqref{e:F1} and \eqref{e:F2} we have 
\[
\begin{array}{lcrr}
G^{+}_1(0)= \sqrt{2\varphi_1''(0)^{-1}}\, ,\ G^{-}_1(0)= -\sqrt{-2\varphi_1''(\omega_1/4)^{-1}}\, ,\\[0.3cm]
G^{+}_2(0)= \sqrt{2\varphi_2''(0)^{-1}}\, ,\ G^{-}_2(0)= -\sqrt{-2\varphi_2''(\omega_2 /4)^{-1}}\, .  
\end{array}
\]  
\end{Remark}
As a corollary we obtain
\begin{Lemma}\hspace{-2mm}{\bf .}\label{Lem:analyticity}
The functions $J_1$, $J_2$, and $J_3$ are analytic in 
$(\ka_1,\ka_2)\in (\nu_3,0)\times (0,\nu_1)$. 
\end{Lemma}
{\em Proof of Lemma \ref{Lem:analyticity}}.  
The function $J_1$ is obviously analytic in that domain. 
Fix $a\in (0,\nu_1)$ and take $0<\delta \ll 1$ such that $\rho_2(z)$ is holomorphic in the disc $\D_{2\delta}(a):=\{|z-a|<2\delta\}\subset \C$. Then write
\[
J_2(\ka_1,\ka_2)=
\int_{0}^{a-\delta}\sqrt{\ka_2-x_2}\,f(x_2,\ka_1)\,dx_2\, + \, \int_{a-\delta}^{\ka_2}\sqrt{\ka_2-x_2}\,f(x_2,\ka_1)\,dx_2\, ,
\]
where
\[
f(x_2,\ka_1)= \frac{1}{\pi}\sqrt{x_2-\ka_1}\,\rho_2(x_2)
\]
is analytic in $(0,\nu_1)\times (\nu_3, 0)$. Then the  first integral defines an analytic function in $(\ka_1,\ka_2)\in (\nu_3,0)\times (a-\delta/2,a+\delta/2)$. Consider now the second one. We expand  $f(x_2,\ka_1)$ in Taylor series with respect to $x_2$ at $x_2=\ka_2$.  Then  integrating with respect to $x_2$ and using Cauchy inequalities for 
$\frac {d^jf}{d x_2^j}(\kappa_2,\kappa_1)$, where  $(\kappa_2,\kappa_1)\in \D_{\delta/2}(a)\times (\nu_3+\delta,-\delta)$,  we obtain  that the second integral defines an analytic function in 
$(\ka_1,\ka_2)\in (\nu_3+\delta,-\delta)\times \D_{\delta/2}(a)$. In the same way we prove that $J_3$ is analytic in 
$(\nu_3,0)\times (0,\nu_1)$. 
\finishproof

In order to obtain suitable formulas for the frequencies of the billiard ball map we proceed as in the Appendix. 
Denote by ${\cal H}(J_1,J_1,J_3)$ the Hamiltonian of the billiard flow expressed in the corresponding 
action-angle  coordinates. Then for any $\ka_1$ and $\ka_2$ such that  
$\nu_3<\ka_1<0$ and $0<\ka_2<\nu_1$, one has
\begin{equation}\label{hamiltonian}
{\cal H}(J_1(\ka_1,\ka_2),J_1(\ka_1,\ka_2),J_3(\ka_1,\ka_2))\equiv 1\,.
\end{equation}
Differentiating \eqref{hamiltonian} with respect to $\ka_1$ and $\ka_2$ we get that
the frequencies $\Omega_1$ and $\Omega_2$ of the billiard ball map satisfy
\[
\left[
\begin{array}{ccc}
\frac{\p J_1}{\p\ka_1}&\frac{\p J_2}{\p\ka_1}&\frac{\p J_3}{\p\ka_1}\\
\frac{\p J_1}{\p\ka_2}&\frac{\p J_2}{\p\ka_2}&\frac{\p J_3}{\p\ka_2}
\end{array}
\right]
\left[
\begin{array}{c}
\Omega_1\\
\Omega_2\\
2\pi
\end{array}
\right]
=0
\]
and therefore (cf. formula \eqref{e:frequencies_relation} in the Appendix)
\begin{equation}\label{om_relation}
\displaystyle\left[
\begin{array}{cc}
\frac{\p J_1}{\p\ka_1}&\frac{\p J_2}{\p\ka_1}\\
\frac{\p J_1}{\p\ka_2}&\frac{\p J_2}{\p\ka_2}
\end{array}
\right]
\left[
\begin{array}{c}
\Omega_1\\
\Omega_2
\end{array}
\right]
=
-2\pi\left[
\begin{array}{c}
\frac{\p J_3}{\p\ka_1}\\
\frac{\p J_3}{\p\ka_2}
\end{array}
\right].
\end{equation}
The latter relation and the formulas for the actions \eqref{action1}-\eqref{action3} lead to
the following formulas for the frequencies
\begin{equation}\label{om}
\displaystyle\Omega_1(\ka_1,\ka_2)=\pi\,\frac{A(\ka_1,\ka_2)}
{D(\ka_1,\ka_2)}\quad  \mbox{and}\quad  \Omega_2(\ka_1,\ka_2)=\pi\,\frac{B(\ka_1,\ka_2)}
{D(\ka_1,\ka_2)}\, ,
\end{equation}
where
\[
\displaystyle A(\ka_1,\ka_2):= 
\int_{\nu_3}^{\ka_1}\int_0^{\ka_2}\!\!
\frac{(x_2-x_3)\rho_2(x_2)\rho_3(x_3)\,dx_2\,dx_3}{\sqrt{(x_2-\ka_1)(\ka_2-x_2)(\ka_1-x_3)(\ka_2-x_3)}}\, ,
\]
\[
B(\ka_1,\ka_2):= \
 \displaystyle  \int_{\nu_3}^{\ka_1}\int_{\nu_1}^{\nu_0}\!\!
\frac{(x_1-x_3)\rho_1(x_1)\rho_3(x_3)\,dx_1\,dx_3}{\sqrt{(x_1-\ka_1)(x_1-\ka_2)(\ka_1-x_3)(\ka_2-x_3)}}
\]
and 
\[
D(\ka_1,\ka_2):= 
 \displaystyle  \int_0^{\ka_2}\int_{\nu_1}^{\nu_0}\!\!
\frac{(x_1-x_2)\rho_1(x_1)\rho_2(x_2)\,dx_1\,dx_2}{\sqrt{(x_1-\ka_1)(x_1-\ka_2)(x_2-\ka_1)(\ka_2-x_2)}}\,.
\]
It follows from Lemma \ref{Lem:analyticity} that $A$, $B$ and $D$ are analytic functions in $(\ka_1,\ka_2)\in (\nu_3,0)\times (0,\nu_1)$. Moreover, $D\neq 0$ in that domain, which implies that $\Omega_1$ and $\Omega_2$ are analytic in $(\ka_1,\ka_2)\in (\nu_3,0)\times (0,\nu_1)$.

Denote by ${\mathcal J}$ the Jacobian of the frequency map $(\ka_1,\ka_2)\mapsto(\Omega_1(\ka_1,\ka_2),\Omega_1(\ka_1,\ka_2))$,
\begin{equation}\label{Jacobian}
{\mathcal J}(\ka_1,\ka_2):=\Big|\frac{\p(\Omega_1,\Omega_2)}{\p(\ka_1,\ka_2)}\Big|=\frac{\pi^2}{D^4}\left|
\begin{array}{cc}
A_{\ka_1}D-AD_{\ka_1}&A_{\ka_2}D-AD_{\ka_2}\\
B_{\ka_1}D-BD_{\ka_1}&B_{\ka_2}D-BD_{\ka_2}\, .
\end{array}
\right|
\end{equation}
Since  ${\mathcal J}(\ka_1,\ka_2)$ is analytic in $(\ka_1,\ka_2)\in (\nu_3,0)\times (0,\nu_1)$, 
either   ${\mathcal J}(\ka_1,\ka_2)\neq 0$  in an open dense subset of $(\nu_3,0)\times (0,\nu_1)$ or 
 \begin{equation}
 \label{e:Z}
 {\mathcal J}(\ka_1,\ka_2)=0\quad  \mbox{for any}\quad (\ka_1,\ka_2)\in (\nu_3,0)\times (0,\nu_1)\, .
 \end{equation}
We are going to compute the limit of ${\mathcal J}(\ka_1,\ka_2)$ as $k_1\to \nu_3+0$. 
To do this we will need the following auxiliary Lemma.
\begin{Lemma}\hspace{-2mm}{\bf .}\label{Lem:estimates}
Let $f(x,\ka)$ be a function on $(a,b)\times(a,b)$ such that $f$ and its partial derivatives $f_x$, $f_\ka$ and
$f_{x\ka}$ exist and are continuous and bounded on $(a,b)\times(a,b)$.
Consider the function  $\displaystyle F(\ka):=\int_a^\ka\frac{f(x,\ka)}{\sqrt{\ka-x}}\,\,dx\,.$
Then 
\begin{itemize}
\item[(a)] $F(\ka)\;=\;2 f(a,\ka)\sqrt{\ka-a}+O(|\ka-a|^{3/2})$
\item[(b)] $\displaystyle F'(\ka)\;=\;\frac{f(a,\ka)}{\sqrt{\ka-a}}+O(\!\sqrt{\ka-a}\;)$
\end{itemize}
where the estimates above are uniform in $\ka\in(a,b)$.
\end{Lemma}
{\em Proof of Lemma \ref{Lem:estimates}.}
An integration by parts leads to
\begin{equation}\label{e:partial_integration}
F(\ka)=-\frac{1}{2}\int_a^{\ka-0}f(x,\ka)\,\,d\sqrt{\ka-x}=2 f(a,\ka)\sqrt{\ka-a}+2\int_a^\ka f_x(x,\ka)\sqrt{\ka-x}\,\,dx
\end{equation}
that together with the boundedness of $f_x$ proves $(a)$.
Differentiating \eqref{e:partial_integration} with respect to $\ka$ and using the boundedness of
$f_x$, $f_\ka$, and $f_{x\ka}$, we prove $(b)$.
\finishproof

The expression for $A(\ka_1,\ka_2)$ can be rewritten in the form
\[
A(\ka_1,\ka_2)=\int_{\nu_3}^{\ka_1}\frac{f(x_3,\ka_1;\ka_2)}{\sqrt{\ka_1-x_3}}\,\,dx_3
\]
where
\[
f(x_3,\ka_1;\ka_2):=\frac{\rho_3(x_3)}{\sqrt{\ka_2-x_3}}\int_0^{\ka_2}\!\!\frac{(x_2-x_3)\rho_2(x_2)\,dx_2}{\sqrt{(x_2-\ka_1)(\ka_2-x_2)}}\,.
\]
For any given $\ka_2\in(0,\nu_1)$ the functions $f(x_3,\ka_1;\ka_2)$ and $\frac{\p f(x_3,\ka_1;\ka_2)}{\p\ka_2}$ satisfy the conditions of Lemma \ref{Lem:estimates}
(with $x\equiv x_3$, $\ka=\ka_1$, $a=\nu_3<0<b<0$) in view of Remark \ref{rem:rho}. Applying the Lemma we get
\begin{equation}\label{e:A_principal_part}
\displaystyle A(\ka_1,\ka_2)=\Big(2\frac{\rho_3(\nu_3)}{\sqrt{\ka_2-\nu_3}}\int_0^{\ka_2}\!\!\frac{\sqrt{x_2-\nu_3}\,\rho_2(x_2)\,dx_2}{\sqrt{\ka_2-x_2}}\Big)\sqrt{\ka_1-\nu_3}+
o(\sqrt{\ka_1-\nu_3})
\end{equation}
\begin{equation}\label{e:A_1_principal_part}
\displaystyle \frac{\p A(\ka_1,\ka_2)}{\p\ka_1}=
\Big(\frac{\rho_3(\nu_3)}{\sqrt{\ka_2-\nu_3}}\int_0^{\ka_2}\!\!\frac{\sqrt{x_2-\nu_3}\,\rho_2(x_2)\,dx_2}{\sqrt{\ka_2-x_2}}\Big)\Big/\sqrt{\ka_1-\nu_3}+
o(1/\sqrt{\ka_1-\nu_3})
\end{equation}
and
\begin{equation}\label{e:A_2_principal_part}
\displaystyle \frac{\p A(\ka_1,\ka_2)}{\p\ka_2}=
2\frac{\p}{\p\ka_2}\Big(\frac{\rho_3(\nu_3)}{\sqrt{\ka_2-\nu_3}}\int_0^{\ka_2}\!\!\frac{\sqrt{x_2-\nu_3}\,\rho_2(x_2)\,dx_2}{\sqrt{\ka_2-x_2}}\Big)\sqrt{\ka_1-\nu_3}+
o(\sqrt{\ka_1-\nu_3})
\end{equation}
as $\ka_1\to\nu_3+0$.
In the same way one obtains
\begin{equation}\label{e:B_principal_part}
\displaystyle B(\ka_1,\ka_2)=\Big(2\frac{\rho_3(\nu_3)}{\sqrt{\ka_2-\nu_3}}\int_{\nu_1}^{\nu_0}\!\!\frac{\sqrt{x_1-\nu_3}\,\rho_1(x_1)\,dx_1}{\sqrt{x_1-\ka_2}}\Big)\sqrt{\ka_1-\nu_3}+
o(\sqrt{\ka_1-\nu_3})
\end{equation}
\begin{equation}\label{e:B_1_principal_part}
\displaystyle \frac{\p B(\ka_1,\ka_2)}{\p\ka_1}=
\Big(\frac{\rho_3(\nu_3)}{\sqrt{\ka_2-\nu_3}}\int_{\nu_1}^{\nu_0}\!\!\frac{\sqrt{x_1-\nu_3}\,\rho_1(x_1)\,dx_1}{\sqrt{x_1-\ka_2}}\Big)\Big/\sqrt{\ka_1-\nu_3}+
o(1/\sqrt{\ka_1-\nu_3})
\end{equation}
and
\begin{equation}\label{e:B_2_principal_part}
\displaystyle \frac{\p B(\ka_1,\ka_2)}{\p\ka_2}=
2\frac{\p}{\p\ka_2}\Big(\frac{\rho_3(\nu_3)}{\sqrt{\ka_2-\nu_3}}\int_{\nu_1}^{\nu_0}\!\!\frac{\sqrt{x_1-\nu_3}\,\rho_1(x_1)\,dx_1}{\sqrt{x_1-\ka_2}}\Big)\sqrt{\ka_1-\nu_3}+
o(\sqrt{\ka_1-\nu_3})
\end{equation}
as $\ka_1\to\nu_3+0$.
Note also that for any $\ka_2\in (0,\nu_1)$, $D(\ka_1,\ka_2)$ is a continuous (even real-analytic) function with respect to
$\ka_1$ on the whole interval $(-\infty,0)$.

Consider the limit $\de(\ka_2):=\lim\limits_{\ka_1\to\nu_3+0}{\pi^{-2}D^3 \mathcal J}(\ka_1,\ka_2)$ for  $\ka_2\in(0,\nu_1)$. 
It follows from \eqref{Jacobian} and \eqref{e:A_principal_part}-\eqref{e:B_2_principal_part} that
\begin{eqnarray}
\pi^{-2}D^3 {\mathcal J}
&=&(A_{\ka_1}B_{\ka_2}-B_{\ka_1}A_{\ka_2})+\frac{D_{\ka_1}}{D}(-AB_{\ka_2}+A_{\ka_2}B)+\frac{D_{\ka_2}}{D}(AB_{\ka_1}-A_{\ka_1}B)\nonumber\\
&=&(A_{\ka_1}B_{\ka_2}-B_{\ka_1}A_{\ka_2})+o(1)\nonumber\\
&=&\displaystyle 2\Big(\frac{\rho_3(\nu_3)}{\sqrt{\ka_2-\nu_3}}\Big)^2
\Big(\int_0^{\ka_2}\frac{\sqrt{x_2-\nu_3}\,\rho_2(x_2)\,dx_2}{\sqrt{\ka_2-x_2}}\Big)^2
\frac{\p}{\p\ka_2}\left(\frac{\int_{\nu_1}^{\nu_0}\frac{\sqrt{x_1-\nu_3}\,\rho_1(x_1)\,dx_1}{\sqrt{x_1-\ka_2}}}
{\int_0^{\ka_2}\frac{\sqrt{x_2-\nu_3}\,\rho_2(x_2)\,dx_2}{\sqrt{\ka_2-x_2}}}\right)+o(1)\nonumber
\end{eqnarray}
as $\ka_1\to\nu_3+0$.
Hence,
\begin{equation}\label{e:De_limit}
\displaystyle \de(\ka_2)=
2\Big(\frac{\rho_3(\nu_3)}{\sqrt{\ka_2-\nu_3}}\Big)^2
\Big(\int_0^{\ka_2}\frac{\sqrt{x_2-\nu_3}\,\rho_2(x_2)\,dx_2}{\sqrt{\ka_2-x_2}}\Big)^2
\frac{\p}{\p\ka_2}\left(\frac{\int_{\nu_1}^{\nu_0}\frac{\sqrt{x_1-\nu_3}\,\rho_1(x_1)\,dx_1}{\sqrt{x_1-\ka_2}}}
{\int_0^{\ka_2}\frac{\sqrt{x_2-\nu_3}\,\rho_2(x_2)\,dx_2}{\sqrt{\ka_2-x_2}}}\right)\,.
\end{equation}
Suppose that \eqref{e:Z} holds. Then 
  $\de(\ka_2) = 0$ for any $\ka_2\in(0,\nu_1)$ and it  follows from
\eqref{e:De_limit} that there is a  constant $C\neq 0$ such that
\begin{equation}\label{e:rotation=const}
\int_{\nu_1}^{\nu_0}\frac{\sqrt{x_1-\nu_3}\,\rho_1(x_1)\,dx_1}{\sqrt{x_1-\ka_2}}=C
\int_0^{\ka_2}\frac{\sqrt{x_2-\nu_3}\,\rho_2(x_2)\,dx_2}{\sqrt{\ka_2-x_2}}
\end{equation}
for any $\ka_2\in(0,\nu_1)$. 
\begin{Lemma}\hspace{-2mm}{\bf .}\label{Lem:rotation_function}
Let $(X,\tg)$ be a Liouville billiard table of classical type. Then the geodesic flow of the restriction
${\tilde l}=\tg|_\Ga$ of the Riemannian metric $\tg$ to the boundary $\Ga$ is completely integrable.
A functionally independent with $\tilde l$ integral of the geodesic flow of $\tilde l$ is given by the restriction $\tI=\tI_2|_\Ga$ of $\tI_2$ to $\Ga$ and the level set $\{\tilde l=1, \tilde I=\ka\}$ is non empty if and only if $\kappa\in [0,\nu_0]$. 
In action-angle coordinates the {\em rotation function} corresponding to the Liouville torus
$\tilde T_\kappa:=\{\tilde l=1, \tilde I=\ka\}$\footnote{This set has two connected components that correspond to two Liouville tori
with the same rotation function.} for $\ka\in(0,\nu_1)$ is
\begin{equation}\label{e:rotation_function}
\rho(\ka)=2\int_{\nu_1}^{\nu_0}\frac{\sqrt{x_1-\nu_3}\,\rho_1(x_1)\,dx_1}{\sqrt{x_1-\ka}}\Big/
\int_0^{\ka}\frac{\sqrt{x_2-\nu_3}\,\rho_2(x_2)\,dx_2}{\sqrt{\ka-x_2}}\,.
\end{equation}
\end{Lemma}
{\em Proof of Lemma \ref{Lem:rotation_function}.}
It follows from the construction of the Liouville billiard tables that the mapping $\si|_{\T^2_N} : \T^2_N\to\Ga$
is a double branched covering of the boundary $\Ga$,  where
\[
\T^2_N=\{(\te_1\,(\mod\om_1),\te_2\,(\mod\om_2),\te_3=N)\}\subset C . 
\]
In the coordinates $\{(\te_1,\te_2)\}$ on $\T^2_N$ we get
the following expressions for the metric $l=(\si|_{\T^2_N})^*{\tilde l}$ and the integral $I=(\si|_{\T^2_N})^*\tI$
\begin{eqnarray}
dl^2&=&(\varphi_1-\varphi_2)\Big((\varphi_1-\nu_3)\,d\te_1^2+(\varphi_2-\nu_3)\,d\te_2^2\Big)\,,\nonumber\\
dI^2&=&(\varphi_1-\varphi_2)\Big(\varphi_2(\varphi_1-\nu_3)\,d\te_1^2+\varphi_1(\varphi_2-\nu_3)\,d\te_2^2\Big)\,.\nonumber
\end{eqnarray}
Applying the Legendre transformation corresponding to $l$ we obtain the following system of equations
for the level set $T_\kappa:=\{l=1, I=\ka\}$,
\begin{eqnarray}
L&=&\frac{1}{\varphi_1-\varphi_2}\Big(\frac{p_1^2}{\varphi_1-\nu_3}+\frac{p_2^2}{\varphi_2-\nu_3}\Big)=1\nonumber\\
I&=&\frac{1}{\varphi_1-\varphi_2}\Big(\varphi_2\frac{p_1^2}{\varphi_1-\nu_3}+
\varphi_1\frac{p_2^2}{\varphi_2-\nu_3}\Big)=\ka\nonumber
\end{eqnarray}
that leads to the following expression of the impulses on $T_\ka$,
\begin{eqnarray}
p_1(\te_1)^2&=&(\varphi_1(\te_1)-\nu_3)\,(\varphi_1(\te_1)-\ka)\ge 0\label{e:p_1}\\[0.3cm]
p_2(\te_2)^2&=&(\varphi_2(\te_2)-\nu_3)\,(\ka-\varphi_2(\te_2))\ge 0\label{e:p_2} . 
\end{eqnarray}
In particular,  $T_\ka \neq \emptyset$ if and only if $\kappa\in [0,\nu_0]$. 
Hence, the projection of  $T_\ka$ into the base $\T^2_N$ is given by the union of the sets
\[
A_\ka':=\{(\te_1,\te_2)\;: \;0\le\te_1\le\om_1,\;-f_2(\ka)\le\te_2\le f_2(\ka)\}
\]
and
\[
A_\ka'':=\{(\te_1,\te_2)\;: \;0\le\te_1\le\om_1,\;-f_2(\ka)+\om_2/2\le\te_2\le f_2(\ka)+\om_2/2\}\,.
\]
As the sets $A_\ka'$ and $A_\ka''$ have the same image under the projection $\si|_{\T^2_N} : \T^2_N\to \Gamma$ we restrict
our attention only to the set $A_\ka'$. It follows from \eqref{e:p_1}-\eqref{e:p_2} that the mapping
$r_+:A_\ka'\to T^*\T^2_N$,
\[
(\te_1,\te_2)\mapsto(\te_1,\te_2;\sqrt{(\varphi_1(\te_1)-\nu_3)\,(\varphi_1(\te_1)-\ka)},
\pm\sqrt{(\varphi_2(\te_2)-\nu_3)\,(\ka-\varphi_2(\te_2))}\,),
\]
parametrizes one of the two connected components of the set $\tilde T_\kappa=\{{\tilde l}=1,\tI=\ka\}\subset T^*X$.
By Liouville-Arnold formula we get the following formulas for the corresponding actions
\begin{eqnarray}
J_1(\ka)&=&\frac{2}{\pi}\int_0^{\om_1}\sqrt{(\varphi_1(\te_1)-\nu_3)\,(\varphi_1(\te_1)-\ka)}\,d\te_1\nonumber\\
&=&\frac{8}{\pi}\int_{\nu_1}^{\nu_0}\sqrt{(x_1-\nu_3)\,(x_1-\ka)}\rho_1(x_1)\,dx_1
\end{eqnarray}
\begin{eqnarray}
J_2(\ka)&=&\frac{2}{\pi}\int_{-f_2(\ka)}^{f_2(\ka)}\sqrt{(\varphi_2(\te_2)-\nu_3)\,(\ka-\varphi_2(\te_2))}\,d\te_2\nonumber\\
&=&\frac{4}{\pi}\int_0^{\ka}\sqrt{(x_2-\nu_3)\,(\ka-x_2)}\rho_2(x_2)\,dx_2
\end{eqnarray}
In the corresponding action-angle coordinates 
the Hamiltonian $L$ becomes $L=L^0(J_1,J_2)$, where $L^0$ is smooth, 
 and the frequency vector $\omega$  of the  invariant torus $\{J_1=c_1,\, J_2=c_2\}$ is 
 $\omega=-\Big(\frac{\partial L^0}{\partial J_1}(c_1,c_2), \frac{\partial L^0}{\partial J_1}(c_1,c_2)\Big)$. 
Then, differentiating the relation
\[
L^0(J_1(\ka),J_2(\ka))\equiv 1
\]
with respect to $\ka\in(0,\nu_1)$ we get \eqref{e:rotation_function}.
\finishproof

\noindent We need the following technical Lemma.
\begin{Lemma}\hspace{-2mm}{\bf .}\label{Lem:limits}
Let $m<0<M$ be real constants, $F_1\in C^1([0,M])$, and $F_2\in C^1([m,0])$. Then 
\begin{equation}\label{e:limit1}
\int_0^M\frac{F_1(\sqrt{t})}{\sqrt{t}\sqrt{t-\alpha}}\,dt=-2F_1(0)\log\sqrt{-\alpha}+O(1)
\end{equation}
and
\begin{equation}\label{e:limit2}
\int_m^\alpha\frac{F_2(\sqrt{-t})}{\sqrt{-t}\sqrt{\alpha-t}}\,dt=2F_2(0)\log\sqrt{-\alpha}+O(1)
\end{equation}
as $\alpha\to 0-0$.
\end{Lemma}
{\em Proof.} We have
\[
\int_0^M\frac{F_1(\sqrt{t})}{\sqrt{t}\sqrt{t-\alpha}}\,dt=2F_1(0)\int_0^{\sqrt{M}}\frac{1}{\sqrt{u^2-\alpha}}\,du
+ O(1) = -2 F_1(0)\log\sqrt{-\alpha}+O(1). 
\]
The proof of \eqref{e:limit2} is similar and we omit it.
\finishproof

\noindent  Lemma \ref{Lem:limits} can be applied to the two integrals in \eqref{e:rotation_function} using  Remark \ref{rem:rho}. In this way we obtain
\[
\int_{\nu_1}^{\nu_0}\frac{\sqrt{x_1-\nu_3}\,\rho_1(x_1)\,dx_1}{\sqrt{x_1-\ka}} = -2G_1^+(0)\sqrt{\nu_1-\nu_3}\log\sqrt{\nu_1-\kappa}+O(1)
\]
and
\[
\int_{0}^{\kappa}\frac{\sqrt{x_2-\nu_3}\,\rho_2(x_2)\,dx_2}{\sqrt{\ka-x_2}} = 2G_2^-(0)\sqrt{\nu_1-\nu_3}\log\sqrt{\nu_1-\kappa}+O(1) .
\]
On the other hand,  Remark \ref{rem:rho} and  assumption $(A_3)$, (2), in Sect. \ref{sec:construction} imply 
\[
G^{-}_2(0)= -\sqrt{-2\varphi_2''(\omega_2 /4)^{-1}}= -\sqrt{2\varphi_1''(0)^{-1}}= -G_1^+(0),
\] 
and by \eqref{e:rotation_function} we obtain 
\begin{equation*}\label{e:the_limit}
\rho(\ka)\to 2\;\;\mbox{as}\;\;\;\ka\to\nu_1-0\,.
\end{equation*}
As by \eqref{e:rotation=const}, $\rho\equiv\const$ we conclude that $\rho\equiv 2$ on the interval $(0,\nu_1)$.
The latter implies that all the geodesics of $\Ga$ lying on a torus $\tilde T_\kappa$ with $\ka\in(0,\nu_1)$ (see Lemma \ref{Lem:rotation_function})
are periodic. Using the analyticity of the billiard table and considering the Poincar\'e map in a tubular neighborhood
of the ``hyperbolic'' level set $\{l=1,I=\nu_1\}$ we obtain  that any geodesics of  $\Ga$ corresponding to some  $\ka\in[\nu_1,\nu_2)$
is periodic as well. As the level sets $\{l=1,I=0\}$ and $\{l=1,I=\nu_0\}$ consists of periodic geodesics
we see that all the geodesics on $\Ga$ are periodic. Hence, the assumption that the Jacobian ${\mathcal J}$ of
the frequency map vanishes in an open subset of   $(\ka_1,\ka_2)\in(\nu_3,0)\times(0,\nu_1)$ implies
that all the geodesics of $\Ga$ are periodic. 
The case $(B)$ can be studied by the same argument. \finishproof

\section{Proof of Theorem \ref{Th:Problem_A} and Theorem \ref{Th:main-isospectral}}\label{sec:proofs}
\setcounter{equation}{0} 

In this section we prove Theorem \ref{Th:main-isospectral} and Theorem \ref{Th:Problem_A} formulated in the
introduction. Let $(X,g)$ be a 3-dimensional analytic Liouville billiard table of classical type such that $\Gamma:=\partial X$ admits at least one non closed geodesic.

 We will prove a more general result than Theorem \ref{Th:main-isospectral} which requires only finite smoothness of $K_t$. Namely, 
fix $d>1/2$ and $\ell > 4[2d]+11$, where $[2d]$ is the entire part of  $2d$ and 
$d$ is  the exponent in  (H$_1$). Denote by $C^\ell(\Gamma,\R)$ the corresponding class of Hölder continuous functions.
\begin{Theorem}\hspace{-2mm}{\bf .}\label{Th:isospectral}
 Let $[0,1]\ni t\mapsto K_t$ be a continuous curve in
$C^\ell(\Ga,\R)$  and suppose that it satisfies $(H_1)$ and $(H_2)$, where $\ell$ and $d$ are fixed as above. If $K_0$ and $K_1$ are invariant with respect to the group of symmetries $G$, then $K_0\equiv K_1$. 
\end{Theorem}

\noindent{\em Proof.}
Given $\alpha>0$ and $\tau>2$ we denote by $\Omega^\tau_\alpha$ the set of all frequencies $(\Omega_1,\Omega_2)\in\R^2$ satisfying the Diophantine condition
\[
\mbox{For any}\   (k_1,k_2,k_3)\in\Z^3\, , (k_1,k_2)\neq (0,0)\,  :\quad \big|\Omega_1 k_1 + \Omega_2 k_2 + k_3\big| \ge
\frac{\alpha}{\big(|k_1| + |k_2|\big)^{\tau}}\, .
\]
Note that the set $\Omega^\tau := \cup_{\alpha>0}\Omega^\tau_\alpha$ is of full Lebegues measure in $\R^2$ for any $\tau>2$ fixed (cf. \cite[Proposition 9.9]{Lazut}). 
Then it follows from Theorem \ref{Th:frequencies} that the subset of $U_1\cup U_2$ filled by invariant tori  $\Lambda$ with frequencies in $\Omega^\tau$  is dense in $U_1\cup U_2$. 
Take $0<\tau-2\ll 1$ so that $\ell > ([2d]+1)(\tau + 2) + 7$. Then 
we apply  \cite[Theorem 1.1]{PT3} for any $\Lambda$ in that family.  By  Remark \ref{rem:m=1} we have 
\begin{equation}\label{e:radon_equality1}
{\mathcal R}_{K_0,\mu}(\Lambda)={\mathcal R}_{K_t,\mu}(\Lambda)
\end{equation}
for any $t\in[0,1]$ and for any torus $\Lambda$ with frequency in $\Omega^\tau$,
where $\mu=\langle \pi^+(\xi),n_g\rangle^{-1}$. 
By continuity we obtain \eqref{e:radon_equality1}
for any Liouville torus $\Lambda$ lying in the part  $U_1\cup U_2$ of $B^*\Ga$ corresponding to the boundary cases.
Finally, Theorem \ref{Th:isospectral} follows from \eqref{e:radon_equality1} and Theorem \ref{Th:R-rigid}.
\finishproof

\noindent{\em Proof of Theorem \ref{Th:Problem_A}.}
Let $(X,g)$ be a 3-dimensional analytic Liouville billiard table of classical type and let $\mu=1$ or $\mu=\langle\pi^+(\xi),n_g\rangle^{-1}$.
Assume that $K\in C(\Ga,\R)$ is invariant with respect to the group of symmetries $G=(\Z/2\Z)^3$ of $\Ga$ and let
the mean value of $\mu\cdot K$ on any periodic orbit of the billiard ball map be zero.
It follows from Theorem \ref{Th:frequencies} that the set filled by Liouville tori $\Lambda$ of the billiard ball map with frequency vectors
$\Omega:=(\Omega_1,\Omega_2)\in{\mathbb Q}\times{\mathbb Q}$ is dense in the part of $B^*\Ga$ corresponding to boundary cases.
Let $\Lambda$ be such a {\em rational} torus. In action-angle coordinates, $\Lambda\cong\R^2/\Z^2$.
There exists $N\in\N$ and two relatively prime numbers $p,q\in\Z$ such that $\Omega\equiv\Big(\frac{p}{N},\frac{q}{N}\Big)\,(\mod \Z^2)$.
Hence, there is an affine change of coordinates on $\R^2/\Z^2$ such that $\Omega\equiv(1/N,0)\,(\mod \Z^2)$. Denote,
\[
{\mathcal D}:=\{(x,y)\,:\,0\le x<1/N, 0\le y\le 1\}.
\]
Using the invariance of $\Lambda$ and of the Leray form on $\Lambda$ with respect to $B^N$ we obtain,
\begin{eqnarray}\label{e:zero}
\int_\Lambda(\mu\cdot K)\,\la=\frac{1}{N}\sum_{k=1}^N\int_{\mathcal D}(B^*)^k(\mu\cdot K)\,\la=\int_{\mathcal D}\Big(\frac{1}{N}\sum_{k=1}^N(B^*)^k(\mu\cdot K)\Big)\,\la=0
\end{eqnarray}
as by assumption the mean $\sum_{k=1}^N(B^*)^k(\mu\cdot K)$ vanishes. Using the density of rational tori $\Lambda$ in boundary cases, 
equality \eqref{e:zero}, and Theorem \ref{Th:R-rigid} we see that $K\equiv 0$.
\finishproof

\section{Appendix: Frequencies of integrable billiard tables}\label{sec:appendix}
\setcounter{equation}{0} 

In this appendix we collect the necessary facts used for the computation of the frequency map
in Sect. \ref{sec:frequencies}. Our main task is to derive formula \eqref{e:frequencies} for
the frequencies of the billiard ball map.

Let $(X,g)$, $n=\dim X\ge 2$, be a billiard table with non-empty locally convex boundary $\Gamma$.
Consider the reflection map at the boundary,
\begin{equation}\label{e:reflection}
\rho : TX|_\Gamma\to TX|_\Gamma,\,\,\,\,\,\,\xi\mapsto\xi-2g(\xi,n_g)n_g,
\end{equation}
where $TX|_\Gamma:=\{\xi\in TX\,:\,\pi(\xi)\in\Gamma\}$ is the restriction of the tangent bundle
to $\Gamma$, $\pi : TX\to X$ is the natural projection onto the base, and $n_g$ is the inward unit normal
to the boundary. The restriction $\rho$ is an involution on $TX|_\Gamma$ the set of fixed point of which
coincides with $T\Gamma\subseteq TX|_\Gamma$. Note that $\rho$ preserves the values of
of the Hamiltonian $H_g(\xi):=\frac{1}{2}\,g(\xi,\xi)$ and when restricted to the unit spherical bundle
$S_gX|_\Gamma:=\{\xi\in TX|_\Gamma\,:\,\|\xi\|_g=1\}$ it coincides with the mapping $r : \Sigma\to\Sigma$
considered in Sect.\,\,\ref{sec:set-up} if we identify vectors and covectors with the help of the Legendre
transform, 
\[
FL_g : TX\to T^*X,\,\,\,\,\,\,\,\xi\mapsto g(\xi,\cdot)\,.
\]
More generally, the notions and mappings considered in Sect.\,\ref{sec:set-up} have their
analogs on $TX$ via the Legendre transform.

Denote by $\al_g$ the Liouville 1-form on $TX$ given by $\al_g(v)(\cdot):=g(v,d_v\pi(\cdot))$
where $v\in TX$ and $(\cdot)$ stands for an arbitrary element of $T_v(TX)$. Note that the differential
$\om_g:=d\al_g$ of the 1-form $\al_g$ corresponds to the symplectic form $dp\wedge dx$ on the cotangent
bundle $T^*X$ via the Legendre transform.
\begin{Lemma}\hspace{-2mm}{\bf .}\label{Lem:reflection_map}
The reflection map $\rho : TX|_\Gamma\to TX|_\Gamma$ satisfies the following properties:
\begin{itemize}
\item[$(a)$] the reflection $\rho$ preserves the restriction of the Liouville form $\al_g$ to $TX|_\Gamma$;
\item[$(b)$] the reflection $\rho$ preserves the values of the Hamiltonian $H_g(\xi)=\frac{1}{2}\,g(\xi,\xi)$;
\item[$(c)$] in the case when $(X,\tg)$ is a Liouville billiard table the reflection $\rho$ corresponding to
the Riemannian metric $\tg$ preserves the values of the pairwise commuting integrals $\tI_k$ $(k=1,2)$ of
the billiard flow (cf. Proposition \ref{Prop:Factorization}).
\end{itemize}
\end{Lemma}
\noindent{\em Proof of Lemma \ref{Lem:reflection_map}.}
$(a)$ Let $t\mapsto v(t)$ be a smooth curve in $TX|_\Gamma$ defined in an open neighborhood of $t=0$ such that
$v(0)=v\in T_xX$, $x\in\Gamma$, and ${\dot v}(0)=\Xi\in T_v(TX|_\Gamma)$. One has
\begin{eqnarray}
\al_g(\rho(v))(d_v\rho(\Xi))&=&g(\rho(v),d_{\rho(v)}\pi\circ d_v\rho(\Xi))=
g\Big(\rho(v),\frac{d}{dt}\pi(\rho(v(t))|_{t=0})\Big)\nonumber\\
&=&g\Big(v-2g(v,n_g)n_g,\frac{d}{dt}\pi(v(t))|_{t=0}\Big)=
g\Big(v,\frac{d}{dt}\pi(v(t))|_{t=0}\Big)\nonumber\\
&=&\al_g(v)(\Xi)
\end{eqnarray}
where we have used that $\pi\circ\rho=\pi$ and that $\frac{d}{dt}\pi(v(t))|_{t=0}\in T_x\Gamma$
is orthogonal to $n_g$. This proves statement $(a)$.
The proof of $(b)$ is straightforward and we omit it.
Statement $(c)$ was established in the proof of Proposition \ref{Prop:Factorization}.
\finishproof

Now we will describe a special variant of the symplectic gluing procedure introduced by Lazutkin
in \cite[\S\,4]{Lazut}. The main idea is to identify parts of the boundary $\partial(TX)$ of the configuration space
$TX$ of the billiard flow in order to ``eliminate" the reflections and obtain a new ``glued'' configuration space
together with a smooth billiard flow on it. Note that the glued configurations space becomes a smooth symplectic manifold
so that the billiard flow is a smooth Hamiltonian system on it. Divide the boundary of $TX$ into three parts
\[
\partial(TX)=TX|_\Gamma=T^-X|_\Gamma\sqcup T^+X|_\Gamma\sqcup T\Gamma
\]
where $T^\pm X|_\Gamma:=\{\xi\in TX_\Gamma\,:\,\pm g(\xi,n_g)>0\}$ and $T\Gamma$ is assumed naturally embedded
into $TX|_\Gamma$.
Note that
\begin{equation}\label{e:reflection+-}
\rho|_{T^-X} : T^-X|_\Gamma\to T^+X|_\Gamma
\end{equation}
is a diffeomorphism and the elements of $T\Gamma\subseteq TX|_\Gamma$ are fixed points of $\rho$.
Now, using \eqref{e:reflection+-} we identify the points $\xi^-\in T^-X$ and $\rho(\xi^-)\in T^+X$ of the boundary
of $TX\setminus T\Gamma$ and obtain a new {\em glued} space $\widetilde{TX}^\rho$ that we supply with
the factor topology so that the projection $\pi_\rho : TX\setminus T\Gamma\to\widetilde{TX}^\rho$,
\[
TX\setminus T\Gamma\ni\xi\mapsto\left\{
\begin{array}{l}
\xi\,\,\,\,\mbox{if}\,\,\,\,\xi\notin\partial(TX)\\
\{\xi,\rho(\xi)\}\,\,\,\,\mbox{if}\,\,\,\,\xi\in T^-X|_\Gamma\\
\{\rho^{-1}(\xi),\xi\}\,\,\,\,\mbox{if}\,\,\,\,\xi\in T^+X|_\Gamma
\end{array}
\right.,
\]
is continuous. 
\begin{Def}\hspace{-2mm}{\bf .}\label{Def:integrability_billiard_flow} 
The billiard flow of $(X,g)$ is called {\em completely integrable} if there exist $n$ functionally
independent integrals $Q_1,...,Q_n\equiv H_g\in C^\infty(TX,\R)$ of the billiard flow such that
$\forall\, 1\le k,l\le n$, $\{Q_k,Q_l\}=0$, and $\forall\, 1\le k\le n$ $\forall\,\xi\in TX|_\Gamma$, 
$Q_k(\rho(\xi))=Q_k(\xi)$.
\end{Def}
Assume that the billiard flow on $TX$ is {\em completely integrable}.
Denote by $X_g$ the Hamiltonian vector field on $TX$ with Hamiltonian $H_g$.
The following Proposition follows from Lemma \ref{Lem:reflection_map} $(a)$, $(b)$, and
is a special case of the symplectic gluing developed in \cite[\S\,4]{Lazut}.
\begin{Prop}\hspace{-2mm}{\bf .}\label{Prop:gluing}
There exists a smooth differentiable structure on $\widetilde{TX}^\rho$, a symplectic form ${\tilde\om}_g$
on $\widetilde{TX}^\rho$, and functions $\tQ_k\in C^\infty(\widetilde{TX}^\rho,\R)$ $(1\le k\le n)$,
such that the projection
\begin{equation}\label{e:factorization_map}
\pi_\rho : TX\setminus T\Gamma\to\widetilde{TX}^\rho
\end{equation}
is smooth, $\pi_\rho^*({\tilde\om}_g)=\om_g$, and $\pi_\rho^*(\tQ_k)=Q_k$ for any $1\le k\le n$.
In particular, the Hamiltonian vector field ${\tilde X}_g$ corresponding to ${\tilde H}_g:=\tQ_n$ is
completely integrable in $\widetilde{TX}^\rho$ and $(\pi_\rho)_*(X_g)={\tilde X}_g$.
\end{Prop}
Denote,
\[
{\mathcal T}:=\pi_\rho(T^\pm X|_\Gamma)\subset\widetilde{TX}^\rho\,.
\]
Note that ${\mathcal T}$ is a disjoint union of connected non-intersecting embedded hypersurfaces in
$\widetilde{TX}^\rho$ that are transversal to the Hamiltonian vector field ${\tilde X}_g$.

Let $c=(c_1,...,c_n)$ be a regular value of the ``momentum'' map
\[
M : \widetilde{TX}^\rho\to\R^n,\,\,\,\,\,\,\,\xi\mapsto(\tQ_1(\xi),...,\tQ_n(\xi))
\] 
and let $\tilde T_c$ be a connected component of the level set $M^{-1}(c)$. The compactness of $X$ implies that
$\tilde T_c$ is compact. By the Liouville-Arnold theorem $\tilde T_c$ is diffeomorphic to the $n$ dimensional
torus $\T^n$ and one can introduce {\em action-angle} coordinates in a tubular neighborhood of $\tilde T_c$ in
$\widetilde{TX}^\rho$ (see \cite{Arn}). 
Assume that 
\[
{\tilde T_c}\cap{\mathcal T}\ne\emptyset
\]
is a non-empty compact set. In this case we will call ${\tilde T_c}$ {\em glued} Liouville torus.
As ${\tilde X}_g$ is tangent to ${\tilde T_c}$ and transversal to ${\mathcal T}$ the submanifolds ${\tilde T_c}$ and
${\mathcal T}$ intersect transversally. Hence, ${\tilde T_c}\cap{\mathcal T}$ is a disjoint union of finitely many
compact embedded submanifolds in ${\tilde T_c}$. Denote by $m\ge 1$ the number of the connected components of
${\tilde T_c}\cap{\mathcal T}$. The proof of the following Lemma is straightforward and we omit it.
\begin{Lemma}\hspace{-2mm}{\bf .}\label{Lem:glued_torus}
\begin{itemize}
\item[$(a)$] The connected components of ${\tilde T_c}\cap{\mathcal T}$ are diffeomorphic to $\T^{n-1}$;
\item[$(b)$] The closure of any of the connected components of ${\tilde T_c}\setminus{\mathcal T}$ is diffeomorphic to
$[0,1]\times\T^{n-1}$. The $n$-torus ${\tilde T_c}$ is obtained by a ``cyclic" gluing together of all $m\ge 1$ copies of
$[0,1]\times\T^{n-1}$ along their boundaries;
\item[$(c)$] Let ${\mathcal S}_c$ be a connected component of ${\tilde T_c}\cap{\mathcal T}$ and let
$\Lambda_c=p_+(\pi_\rho^{-1}({\mathcal S}_c))$ where $p_+ : TX|_\Ga\to T\Ga$ denotes the orthogonal projection
$\xi\mapsto\xi-g(n_g,\xi)\,n_g$ onto $T\Ga$.
Then the number $m\ge 1$ of the connected components of ${\tilde T_c}\cap{\mathcal T}$ is the minimal power of
the billiard ball map $B$ that leaves $\Lambda_c$ invariant,
i.e., $B^m(\Lambda_c)=\Lambda_c$.\footnote{Note that we identify vectors and covectors via the Riemannian metric $g$.}
\end{itemize}
\end{Lemma}
Choose a component ${\mathcal S}_c$ of ${\tilde T_c}\cap{\mathcal T}$ and a basis of cycles
${\tilde\ga}_1,...,{\tilde\ga}_{n-1}$ of its homology group as well as a {\em transversal} cycle
${\tilde\ga}_n$ in ${\tilde T_c}$ so that ${\tilde\ga}_1,...,{\tilde\ga}_n$ is a basis of the homology group of
${\tilde T_c}$. Let $\{{\tilde J}_1,...,{\tilde J}_n;{\tilde\te}_1(\mod2\pi),...,{\tilde\te}_n(\mod2\pi)\}$ be
action-angle coordinates in a tubular neighborhood of the glued Liouville torus ${\tilde T_c}$ that corresponds
the the cycles ${\tilde\ga}_1,...,{\tilde\ga}_n$, i.e., $\forall\,1\le k\le n$,
\[
{\tilde J}_k=\frac{1}{2\pi}\int_{\tilde\ga_k}{\tilde\al}_g\,,
\]
where ${\tilde\al}_g:=(\pi_\rho)_*(\al_g)$ is the push-forward of the Liouville form $\al_g$ onto
$\widetilde{TX}^\rho$. In the action-angle coordinates,
\[
{\tilde X}_g=\eta_1(\tilde J)\frac{\p}{\p{\tilde\te}_1}+...+\eta_n(\tilde J)\frac{\p}{\p{\tilde\te}_n}\,,
\]
where 
\begin{equation}
\label{e:the-eta}
\eta_k({\tilde J}_1,...,{\tilde J}_n):=
\frac{\p\widetilde{\mathcal H}_g}{\p{\tilde J}_k}({\tilde J}_1,...,{\tilde J}_n)\, ,\ 1\le k\le n.
\end{equation}
It follows from the choice of the cycles ${\tilde\ga}_1,...,{\tilde\ga}_n$ that ${\mathcal S}_c$ is a section of the bundle,
\[
\T^{n-1}\times\T\to\T^{n-1},\,\,\,\,\,\,\,\,\,\,
({\tilde\te}_1,...,{\tilde\te}_n)\mapsto({\tilde\te}_1,...,{\tilde\te}_{n-1})\,.
\]
As ${\mathcal S}_c$ is transversal to ${\tilde X}_g$ one concludes that $\eta_n(\tilde J)\ne 0$.
It follows from our construction that the billiard ball map $B^m$ is conjugated to the following
diffeomorphism of the $n-1$-dimensional torus $\{({\tilde\te}_1(\mod2\pi),...,{\tilde\te}_{n-1}(\mod2\pi))\}$,
\[
{\tilde\te}_k\mapsto{\tilde\te}_k+2\pi\frac{\eta_k(\tilde J)}{\eta_n(\tilde J)},\,\,\,\,1\le k\le n-1.
\]
Parameterizing the glued Liouville tori with fixed energy $\{{\tilde H}_g=1\}$ by the values of the integrals
$\tQ=(\tQ_1,...,\tQ_{n-1})$ we obtain the following mapping for the frequencies of $B^m$,
\begin{equation}\label{e:frequency_map}
(\tQ_1,...,\tQ_{n-1})\mapsto(\Omega_k(\tQ_1,...,\tQ_{n-1}))_{1\le k\le n-1}:=
2\pi\,\left(\frac{\eta_k(J(\tQ_1,...,\tQ_{n-1},1))}{\eta_n(J(\tQ_1,...,\tQ_{n-1},1))}\right)_{1\le k\le n-1}
\end{equation}
where $\eta_k({\tilde J}_1,...,{\tilde J}_n)$ is defined by \eqref{e:the-eta}.
Finally, by partial differentiation of the identity,
\[
\widetilde{\mathcal H}_g({\tilde J}(\tQ_1,...,\tQ_{n-1},1))\equiv 1,
\]
one gets that the frequency vector $\Omega:=(\Omega_1,...,\Omega_{n-1})^T$ satisfies the
linear relation
\begin{equation}\label{e:frequencies_relation}
{\bf A}\Omega=-2\pi\,{\bf b}
\end{equation}
where ${\bf A}(Q):=\Big(\frac{\p J_l}{\p Q_k}(Q,1)\Big)_{1\le k,l\le n-1}$,
${\bf b}(Q):=(\frac{\p J_n}{\p Q_1}(Q,1),...,\frac{\p J_n}{\p Q_{n-1}}(Q,1))^T$, and $\forall 1\le k\le n$,
\begin{equation}\label{actions}
J_k:=\frac{1}{2\pi}\int_{\ga_k}\al_g\,,
\end{equation}
where $\ga_k$ is the connected component of $\pi_\rho^{-1}({\tilde\ga}_k)$ lying in $T^+X|_\Gamma$ and
$Q:=(Q_1,...,Q_{n-1})$. The functions $J_k$ $(1\le k\le n)$ will be called {\em generalized}
actions of the billiard flow. Using that $\eta_n\ne 0$ one can prove that ${\bf A}(Q)$ is non-degenerate. Hence,
\begin{equation}\label{e:frequencies}
\Omega(Q)=-2\pi\,{\bf A}(Q)^{-1}\,{\bf b}(Q),
\end{equation}
where $Q_1,...,Q_{n-1}$ are the integrals of the billiard flow in a tubular neighborhood of
the invariant set $\pi_\rho^{-1}({\tilde T}_c)$ of the billiard flow.

\vspace{0.5cm} 
\noindent

\vspace{0.5cm} 
\noindent 
G. P.: 
Universit\'e de Nantes,  \\
Laboratoire de mathématiques Jean Leray,
CNRS: UMR 6629,\\
2, rue de la Houssini\`ere,\\
 BP 92208,  44072 Nantes 
Cedex 03, France \\

\vspace{0.5cm} 
\noindent P. T.: 
Northeastern University,\\
 Department of Mathematics,\\
  360 Huntington Avenue,
Boston, MA, 02115


\begin{thebibliography}{99} 

\bibitem{Adams} {\sc J. Adams}: {\em Expression of the product of any two Legendre's
coefficients by means of Legendre's coefficients}, Proc. R. Soc. Lond, $\bf 27$(1878),
63-71

\bibitem{Arn} {\sc V. Arnold}: {\em Mathematical methods of classical mechanics},
Springer-Verlag, NY, 1989

\bibitem{Besse} {\sc A. Besse}: {\em Manifolds all of whose geodesics are closed},
Springer-Verlag, Berlin-New York, 1978

\bibitem{GM} {\sc V. Guillemin and R. Melrose}: 
{\em An inverse spectral result for 
elliptical regions in ${\mathbb R}^2$}, 
Advances in Mathematics, $\bf 32$(1979), 128-148 
 
\bibitem{GM1} {\sc V. Guillemin and R. Melrose}: 
{\em The Poisson summation formula for manifolds with boundary}, 
Advances in Mathematics, $\bf 32$(1979), 204-232 
 
\bibitem{Kiyo} {\sc K. Kiyohara}: 
{\em Two classes of Riemannian manifolds whose geodesic flows are integrable}, 
Memoirs of the AMS, $\bf 130$(1997), Number 619 

\bibitem{Knorrer} {\sc H. Kn\"orrer}:
{\em Singular fibers of the momentum mapping for integrable Hamiltonian systems}, 
J. Reine Angew. Math., $\bf 355$(1985), 67-107
  
\bibitem{Lazut} {\sc V. Lazutkin}:
{\em KAM theory and semiclassical approximations to eigenfunctions},
Springer-Verlag, Berlin, 1993  

\bibitem{Mos-Ves} {\sc J. Moser and A. Veselov}: {\em Discrete  versions of 
some integrable systems and factorization of matrix polynomials}, 
Comm. Math. Phys., $\bf 139$(1991), 217-243 
   
   
\bibitem{PT1} {\sc G. Popov and P. Topalov}:
{\em Liouville billiard tables and an inverse spectral result},
Ergod. Th. \& Dynam. Sys., $\bf 23$(2003), 225-248

\bibitem{PT2} {\sc G. Popov and P. Topalov}:
{\it   Discrete analog of the projective equivalence and integrable billiard tables}, 
Ergod. Th. \& Dymam. Sys., $\bf 28$ (2008),  1657-1684.

\bibitem{PT3} {\sc G. Popov and P. Topalov}:
{\em Invariants of isospectral deformations and spectral rigidity}, preprint 2009. 

\bibitem{TopCrit} {\sc P. Topalov}: 
{\it Integrability criterion of geodesical 
equivalence. Hierarchies}, Acta Appl. Math., 
$\bf 59$(3)(1999), 271-298 
 
 
 
\bibitem{Perelomov} {\sc A. Perelomov}:
{\em Integrable Systems of Classical Mechanics and Lie Algebras},
Birkh{\"a}user-Verlag, Basel, 1990

 
 
\bibitem{Tabach} {\sc S. Tabachnikov}: 
{\em Billiards}, Panoramas et Syntheses, 
Societe Mathematique de France, 1995 
 
\bibitem{WW} {\sc E. Whittaker and G. Watson}:
{\em A course of modern analysis}, Cambridge University Press, 1927

\end{thebibliography}
\end{document}